\documentclass[preprint,12pt]{elsarticle}

\usepackage{url}
\usepackage{graphicx}

\usepackage{algorithm}
\usepackage{algorithmic}
\usepackage{amssymb}
\usepackage{amsmath}
\usepackage{epstopdf}
\usepackage{xcolor}

\usepackage{subeqnarray}

\newtheorem{remark}{Remark}[section]

\definecolor{grey}{rgb}{.7,.6,.5}

\definecolor{darkgreen}{rgb}{0.05, 0.5, 0.06}

\definecolor{darkred}{rgb}{0.5, 0.05, 0.06}

\begin{document}

\begin{frontmatter}

\title{Digital Twins: McKean--Pontryagin Control for Partially Observed Physical Twins}
\author[label1]{Manfred Opper} 
\ead{manfred.opper@tu-berlin.de}
\author[label2]{Sebastian Reich}
\ead{sebastian.reich@uni-potsdam.de}

\affiliation[label1]{
            organization={Institut f\"ur Softwaretechnik und Theoretische Informatik, Technische Universit\"at Berlin},
            %addressline={}, 
            city={Berlin},
            postcode={D-10587}, 
            %state={},
            country={Germany}}
\affiliation[label2]{organization={Institut f\"ur Mathematik, Universit\"at Potsdam},
            %addressline={}, 
            city={Potsdam},
            postcode={D-14476}, 
            %state={},
            country={Germany}}

%\maketitle

\begin{abstract} 
Optimal control for fully observed diffusion processes is well established and has led to numerous numerical implementations based on, for example, Bellman's principle, model free reinforcement learning, Pontryagin's maximum principle, and model predictive control. In contrast, much fewer algorithms are available for optimal control of partially observed processes. However, this scenario is central to the digital twin paradigm, where a physical twin is partially observed and control laws are derived based on a digital twin. In this paper, we contribute to this challenge by combining data assimilation in the form of the ensemble Kalman filter with the recently proposed McKean--Pontryagin approach to stochastic optimal control. We derive forward evolving mean-field evolution equations for states and co-states which simultaneously allow for an online assimilation of data as well as an online computation of control laws. The proposed methodology is therefore perfectly suited for real time applications of digital twins. We present numerical results for controlled Lorenz-63 and Lorenz-96 systems as well as an inverted pendulum.
\end{abstract}


\end{frontmatter}

%%%%%%%%%%%%%%%%%%%%%%%%%%%%%%%%%%%%%%%%%%%%%
%
\section{Introduction}
%
%%%%%%%%%%%%%%%%%%%%%%%%%%%%%%%%%%%%%%%%%%%%%%

In this paper, we consider digital twins in the form of stochastic processes, which are used to derive appropriate control laws applicable to its physical twin. The physical twin in turn is partially observed subject to measurement errors. To formalize the problem mathematically, we introduce a physical process (physical twin) as a controlled stochastic differential equation
of the form
\begin{equation} \label{eq:FSDE}
\dot{X}_t = b(X_t) + G(X_t) U_t  + \Sigma^{1/2} \dot{B}_t.
\end{equation}
The goal is to find controls $U_t$ that minimizes an infinite horizon discounted cost function
\begin{equation} \label{eq:cost}
J(U) = \mathbb{E} \left[ \int_0^\infty e^{-\gamma t}\left(c(X_t) + \frac{1}{2} \|U_t\|^2 \right) {\rm d}t  \right]
\end{equation}
Here, $B_t$ denotes $d_x$-dimensional Brownian motion, $\Sigma \in \mathbb{R}^{d_x\times d_x}$ the symmetric positive definite diffusion matrix, $G(x) \in \mathbb{R}^{d_x\times d_u}$ the possibly position-dependent control matrix, $c(x)\ge 0$ the running cost, and $\gamma >0$ the discount factor. Expectations are taken with respect to the law $\pi_t$ of the process $X_t$. 

Infinite horizon cost functions are typically considered in model-free \cite{Meyn} as well as model-based reinforcement learning \cite{Model-based_reinforcement} under the assumption of fully observed physical twins. The purpose of this contribution is instead to develop a mathematical and computational framework to approximate controls $U_t$ that depend on the law $\pi_t$ of $X_t$. Such a scenario arises, for example, from partially observed physical twins for which their states, denoted by $X_t^\dagger$, are unknown and can only be estimated via indirect measurements
\begin{equation} \label{eq:obs}
Y_{t_n}^\dagger = h(X^\dagger_{t_n}) + \Xi_{t_n}^\dagger
\end{equation}
in intervals $\Delta \tau>0$ at discrete times $t_{n} = n\Delta \tau$, $n\ge 1$, subject to Gaussian measurement noise $\Xi_{t_n}^\dagger \sim {\rm N}(0,R)$ with covariance matrix $R \in \mathbb{R}^{d_y\times d_y}$. Here $h(x) \in \mathbb{R}^{d_y}$ denotes the forward map that links the unknown states $X_t^\dagger \in \mathbb{R}^{d_x}$ of the physical twin to the observed data $Y^\dagger_{t_n}\in \mathbb{R}^{d_y}$. 

Upon combining these partial and noisy observations with the stochastic differential equation model (\ref{eq:FSDE}) in a process called data assimilation \cite{reich2015probabilistic,law2015data}, our knowledge about $X_t^\dagger$ is then captured by the conditional distributions $\pi_t(x|\{Y^\dagger_{t_n}\}_{t_n\le t})$. Thus, the desired controls $U_t$ can only depend on these distributions, which are also called belief states in the literature on partially observed Markov decision processes \cite{ASTROM1965,Bensoussan92,Handel07,Meyn}, and not on the states $X_t^\dagger$ of the physical twin themselves. 

The mathematical formulation of a partially observed Markov decision processes in continuous time leads either to a Hamilton--Jacobi--Bellman equation in the value functional $v_t(\pi)$, where $\pi$ is a probability distribution \citep{Nisio15}; or, alternatively, to a pair of forward and backward  stochastic partial differential equations \citep{Bensoussan92}. Scalable algorithms for solving the Hamilton--Jacobi--Bellman equation in $v_t(\pi)$ or forward-backward stochastic partial differential equations are currently unavailable. 

Indeed, most current algorithms for partially observed Markov decision processes assume either discrete state and action spaces \citep{Silver10} or assume variants of the separation/ certainty equivalence principle \cite{Handel07,WW81}; {\it i.e.}, find the optimal control for fully observed processes, denoted by $u_\ast (x)$, and then use 
\begin{equation}
    U_t = \int_{\mathbb{R}^{d_x}} u_\ast(x)\,\pi_t(x|\{Y_{t_n\le t}^\dagger\}){\rm d}x
\end{equation}
as a control in the partially observed setting \citep{QMDP,RPPC08}. In an alternative approach, future data are accounted for via local Gaussian approximations giving rise to $\pi_t$-dependent linear control laws \citep{9811560,BSP24}. Very recently, the optimal control-as-inference approach \citep{Kappen05,KVO12} has been extended to partially observed Markov decision processes using sequential Monte Carlo techniques \citep{AIS25}. The ensemble Kalman filter has been combined with model predictive control for partially observed processes in \cite{KK24,R25a}.

In this paper, we instead combine the ensemble Kalman filter \cite{Evensenetal2022} for assimilating data with the McKean--Pontryagin approach to infinite horizon optimal control
\cite{R25b} in order to derive an interacting particle formulation in states $X_t$ and co-states $P_t$, which allows for an online approximation of the desired control law $U_t$ for partially observed physical twins. In other words, we propose a digital twin in the form of evolution equations in terms of $M>1$ interacting particles $(X_t^{(i)},P_t^{(i)}) \in \mathbb{R}^{2d_x}$, $i=1,\ldots,M$. The physical and digital twins interact through the data $Y_{t_n}^\dagger$ of the physical twin and the controls $U_t$ derived from the digital twin and applied to its physical twin.

The remainder of this contribution is structured as follows. The mathematical background required on the ensemble Kalman filter \cite{Evensenetal2022,CRS22} for discrete- and continuous-time observations is summarized in Section \ref{sec:background}. In the same section, we also summarize the McKean--Pontryagin formulation of infinite horizon stochastic optimal control from \cite{R25b}. Both the ensemble Kalman filter and the McKean--Pontryagin formulation are combined in Section \ref{sec:McKean--Pontryagin} and deliver our novel approach to the design of digital twins for online control of partially observed physical twins. Numerical implementation aspects are discussed in Section \ref{sec:implementation}. The controlled Lorenz-63 proposed in \cite{KK24}, a controlled Lorenz-96 model \cite{lorenz96}, and an inverted pendulum \cite{Meyn} are used in Section \ref{sec:example} to illustrate the performance of the proposed methodology. In particular, we study the impact of limiting the absolute value of the control on the ability of the control to force solutions of the Lorenz-63 system to ${\rm x}$-values larger than zero. The inverted pendulum leads to a highly nonlinear control problem while the chaotic Lorenz-96 model is used to demonstrate that the proposed methodology is applicable to higher-dimensional control problems. The paper concludes with some remarks on open problems and related work.

%%%%%%%%%%%%%%%%%%%%%%%%%%%%%%%%%%%%%%
%
\section{Mathematical Background} \label{sec:background}
%
%%%%%%%%%%%%%%%%%%%%%%%%%%%%%%%%%%%%%

In this paper, we assume that the physical twin is represented by (\ref{eq:FSDE}) for given controls $U_t$. It is important to keep in mind that the initial state
$X_0 = X_0^\dagger$ of the physical twin (\ref{eq:FSDE}) is unknown, as well as the specific realization of the Brownian motion $B_t = B_t^\dagger$, which leads to evolving states $X_t = X_t^\dagger$ in time $t\ge 0$ of the physical twin. We use $\dagger$ to distinguish this physical realization from realizations obtained from a digital twin. 

In the following two subsections, we first describe how the ensemble Kalman filter can be used to adjust a digital twin to the available data (\ref{eq:obs}). We then summarize the standard stochastic optimal control problem under the assumption that the physical twin can be fully observed. Both aspects will later be combined to address the control problem for partially observed physical twins.

%%%%%%%%%%%%%%%%%%%%%%%%%%%%%%%%%%%%%
\subsection{Data Assimilation} \label{sec:EnKF}
%%%%%%%%%%%%%%%%%%%%%%%%%%%%%%%%%%%%%

In this paper, we rely on the ensemble Kalman filter \cite{Evensenetal2022,CRS22} to assimilate data into the controlled stochastic differential equation (\ref{eq:FSDE}) for a given control $U_t$. We first consider the case of continuous time data and then return to the discrete-time case (\ref{eq:obs}) subsequently.

\subsubsection{Ensemble Kalman--Bucy Filter}

For ease of presentation, we first consider time-continuous data $Y_t^\dagger$ that satisfy the stochastic differential equation
\begin{equation} \label{eq:cobs}
    \dot{Y}_t^\dagger = h(X_t^\dagger) + R^{1/2}\dot{W}_t
\end{equation}
instead of discrete-time observations (\ref{eq:obs}). Here $W_t \in \mathbb{R}^{d_y}$ denotes Brownian noise independent of $B_t$ and we set $Y_0^\dagger = 0$. In practice, data $Y_t^\dagger$ are observed in small time intervals $\tau>0$ giving rise to
\begin{equation}
\Delta Y_{k\tau}^\dagger := Y_{(k+1)\tau}^\dagger - Y_{k\tau}^\dagger \approx h(X_{k\tau})\tau + (\tau R)^{1/2} \Theta_k
\end{equation}
for $k\ge 1$ and $\Theta_k \sim {\rm N}(0,I)$. Linear interpolation leads to continuous-time and piece-wise differentiable data
\begin{equation} \label{eq:obsc2}
    Y_t^\dagger = \sum_{k=1}^{k_\ast} \Delta Y_{k\tau}^\dagger + \frac{t-k_\ast \tau}{\tau}\Delta Y_{k_\ast \tau}^\dagger
\end{equation}
for $t \in [k_\ast \tau,(k_\ast+1)\tau)$, which is used in our data assimilation and optimal control formulations instead of (\ref{eq:cobs}) since it allows us to treat $\dot{Y}_t^\dagger$ as a regular time derivative.

Following the approach of \cite{bergemann2012ensemble}, we use the following mean-field formulation of the ensemble Kalman--Bucy filter:
\begin{equation}\label{eq:EnKF}
\dot{X}_t = b(X_t) + G(X_t) U_t  + \Sigma^{1/2} \dot{B}_t +  C^{xh}_t R^{-1} \left(\dot{Y}_t^\dagger - \frac{1}{2} \left(
h(X_t) + m^h_t \right)  \right) .
\end{equation}
Here, $C^{xh}_t$ denotes the covariance matrix between $X_t$ and $h(X_t)$ and $m^h_t$ the mean of $h(X_t)$. Provided the controls $U_t$ are known, the mean-field formulation (\ref{eq:EnKF}) provides a digital twin, which can be implemented as an interacting particle system \cite{CRS22}. We note that replacing (\ref{eq:cobs}) by (\ref{eq:obsc2}) leads to a modified Kalman--Bucy filter formulation for nonlinear forward maps \cite{CNN2021}. We ignore the additional drift term in this paper since all numerical experiments rely on linear forward maps.

We recall that the ensemble Kalman--Bucy filter and its mean-field formulation (\ref{eq:EnKF}) only provide approximations to the true filtering distributions. We denote the law of $X_t$, as defined by (\ref{eq:EnKF}), by $\pi^{\rm KBF}_t(x)$. Alternatively, sequential Monte Carlo methods \cite{chopin:20} could be used which, however, are less scaleable to high-dimensional diffusion processes \cite{SBBA08}.

\subsubsection{Ensemble Kalman Filter}

We now return to discrete-time observations (\ref{eq:obs}) for which, following \cite{bergemann2010mollified}, the mean-field ensemble Kalman filter formulation becomes
%\begin{subequations} 
\begin{align} \nonumber \label{eq:mollied EnKF}
\dot{X}_t =& \,\,b(X_t) + G(X_t) U_t  + \Sigma^{1/2} \dot{B}_t\\
& + \sum_{n\ge 1} C^{xh}_t R^{-1} \left(Y_{t_n}^\dagger - \frac{1}{2} \left(
h(X_t) + m^h_t \right) \right) \delta_{t_n}(t),
\end{align}
%\end{subequations}
which can subsequently be implemented as an interacting particle system \cite{CRS22}. Here,
$\delta_{t_n}(t)$ stands for the Dirac delta function centered on $t_n$. We note that
\begin{equation}
\dot{X}_t = C^{xh}_t R^{-1} \left(Y_{t_n}^\dagger - \frac{1}{2} \left(
h(X_t) + m^h_t \right) \right) \delta_{t_n}(t)
\end{equation}
leads to an impulse-like change in $X_t$ at observation time $t=t_n$ while $\dot{X}_t= 0$
otherwise. Let us denote the incoming state at $t_n$ by $X^-_{t_n}$ and the resulting state after the assimilation of $Y_{t_n}^\dagger$ by $X^+_{t_n}$. They satisfy the integral relation
\begin{equation} 
    \tilde{X}_\tau = X_{t_n}^- + \int_0^\tau C^{xh}_s R^{-1} \left(Y_{t_n}^\dagger - \frac{1}{2} \left(h(\tilde X_s) + m^h_s \right) \right) {\rm d}s
\end{equation}
and $X^+_{t_{n+1}} = \tilde X_1$. The integration interval is $\tau \in [0,1]$ since we integrate across the Dirac delta function $\delta_{t_n}$ and
the initial condition is $\tilde{X}_0 = X_{t_n}^-$.
Formulated in terms of Bayesian inference, $X_{t_n}^-$ represents the prior, while $X^+_{t_n}$  encodes the posterior given the data $Y_{t_n}^\dagger$ and the forward model (\ref{eq:obs}).

Again we stress that (\ref{eq:mollied EnKF}) only provides an approximation to the true filtering distributions and denote the law of $X_t$, as defined by (\ref{eq:mollied EnKF}), by
$\pi_t^{\rm KF}(x)$.

%%%%%%%%%%%%%%%%%%%%%%%%%%%%%%%%%%%%%
\subsection{Stochastic optimal control} \label{sec:SOC}
%%%%%%%%%%%%%%%%%%%%%%%%%%%%%%%%%%%%%

The optimal control problem with cost function (\ref{eq:cost}) for fully observed states is solved by the stationary solution $v_\ast (x)$ of the (forward) Hamilton--Jacobi--Bellman equation
\begin{equation}  \label{eq:HJB}
\partial_t v_t = -\gamma v_t + \min_u \left\{ (b+Gu)^{\rm T} \nabla_x v_t + \frac{1}{2} \|u\|^2\right\} + \frac{1}{2}  \Sigma : D_x^2 v_t + c 
\end{equation}
in the value function $v_t(x)$ for $t\ge 0$ with initial condition $v_0 \equiv 0$ \cite{Carmona}; {\it i.e.}
\begin{equation}
v_\ast = \lim_{t\to \infty}v_t.
\end{equation}
Here we used the notation $A:B = \sum_{i,j} a_{ij}b_{ij}$ for suitable matrices $A$ and $B$ \cite{Pavliotis2016}. The optimal closed loop control law is provided by
\begin{equation} \label{eq:optimal control law}
u_\ast(x) = - G(x)^{\rm T} \nabla_x v_\ast (x).
\end{equation}

We note that the Bellman optimality principle states that
\begin{equation} \label{eq:Bellman optimality}
    v_\ast(x) = \min_{u}\mathbb{E} \left[ \int_{0}^{\tau} e^{-\gamma t}\left(c(X_t) + \frac{1}{2} \|U_t\|^2 \right) {\rm d}t + e^{-\gamma \tau}
    v_\ast (X_{\tau}) \right]
\end{equation}
for suitable $\tau > 0$ and with expectations taken with respect to solutions of (\ref{eq:FSDE})  for given control $U_t$, $t\in [0,\tau]$, and initial condition $X_{0} = x$. Formulation (\ref{eq:Bellman optimality}) provides the starting point for numerous algorithms developed in the reinforcement learning community for fully observable physical twins \cite{meyn13}.

It is well-known that the optimal control law can also be found from the Pontryagin minimum (maximum) principle \cite{pontryagin}. More precisely, the classical Pontryagin minimum principle for controlled ordinary differential equations has been extended to controlled stochastic differential equations (\ref{eq:FSDE}) and leads to forward--backward stochastic differential equations in the states $X_t \in \mathbb{R}^{d_x}$, co-states $P_t \in \mathbb{R}^{d_x}$, and Lagrange multipliers $V_t \in \mathbb{R}^{d_x\times d_x}$ \cite{Bensoussan,Carmona}. However, it is computationally non-trivial to apply the stochastic Pontryagin principle directly to infinite horizon optimal control problems since the integration interval for the underlying boundary value problem also becomes infinite.

Recently, an alternative, deterministic mean-field reformulation of the classical Pontryagin minimum principle has been provided in \cite{R25b}. The so called McKean--Pontrygin formulation can be extended to infinite horizon optimal control problems leading to a pair of mean-field ordinary differential equations that are both solved forward in time. In the following, we summarize the key formulas from \cite{R25b}.

Upon introducing functions $\psi_t:\mathbb{R}^{d_x} \to \mathbb{R}^{d_x}$ through the relation
\begin{equation} \label{eq:constraint}
    P_t(a) = \psi_t(X_t(a))
\end{equation}
between states $X_t(a)$ and co-states $P_t(a)$ with labels $a\sim \pi_0$, the McKean--Pontryagin formulation \cite{R25b} results in mean-field evolution equations 
\begin{subequations} \label{eq:PMP}
    \begin{align}
    \dot{X}_t(a) &= \nabla_p H(X_t(a),P_t(a),U_t(a),\beta_t(a);\psi_t) ,\\
    \epsilon \dot{P}_t(a) &= -\gamma P_t(a) + \nabla_x H(X_t(a),P_t(a),U_t(a),\beta_t(a);\psi_t) \nonumber \\
    & \quad \,\,+ \,(1+\epsilon) D_x \psi_t(X_t(a))\dot{X}_t(a),\\
    0 &= \nabla_u H(X_t(a),P_t(a),U_t(a),\beta_t(a);\psi_t),
    \end{align}
\end{subequations}
for a given function $\beta_t(a)$, parameter $\epsilon >0$, and Hamiltonian density
$H(\cdot\,;\psi):\mathbb{R}^{4d_x} \to \mathbb{R}$ defined by
%\begin{subequations} 
\begin{align} \nonumber \label{eq:Hamiltonian density}
H(x,p,u,\beta;\psi) :=&\,  p^{\rm T} (b(x)+G(x)u)   + c(x) +\frac{1}{2} \|u\|^2 \\
& \,+ \beta^{\rm T} \left(\psi(x)-p\right) + \frac{1}{2}\nabla_x \cdot\left(\Sigma \psi(x)\right).
\end{align}
%\end{subequations}
The initial conditions are provided by $X_0(a) = a\sim \pi_0$ and $P_0(a) = 0$. Furthermore, the constraint (\ref{eq:constraint}) follows from
\begin{equation}
0 = \nabla_\beta H(X_t(a),P_t(a),U_t(a),\beta_t(a);\psi_t).
\end{equation}

We note that (\ref{eq:PMP}) was derived in \cite{R25b} for $\epsilon = 1$. However, with sufficient regularity of the stationary value function $v_\ast (x)$ of (\ref{eq:HJB}), the time evolution of the function $\psi_t(x)$ resulting from (\ref{eq:PMP}); {\it i.e.},
\begin{equation}\label{eq:gradient hjb}
    \epsilon \,\partial_t \psi_t = -\gamma \psi_t + \min_u \nabla_x \left(
    (b+Gu)^{\rm T} \psi_t + \frac{1}{2}\|u\|^2 +  c + \frac{1}{2}
    \nabla_x \cdot (\Sigma \psi_t)\right)
\end{equation}
implies
\begin{equation} \label{eq:identity}
\lim_{t\to \infty} \psi_t(x) = \nabla_x v_\ast(x)
\end{equation}
for any choice of $\epsilon >0$.  Hence, the optimal control law (\ref{eq:optimal control law}) can be expressed as
\begin{equation}
    u_\ast(X_t(a)) =  -G(X_t(a))^{\rm T} P_t(a) = -G(X_t(a))^{\rm T}\psi_\ast(X_t(a))
\end{equation}
for $t\to \infty$ with
\begin{equation}
\psi_\ast = \lim_{t\to \infty}\psi_t.
\end{equation}
Furthermore, (\ref{eq:PMP}) can be approximated by
\begin{subequations} \label{eq:PMPs}
    \begin{align}
    \dot{X}_t &= \nabla_p H(X_t,P_t,U_t,\beta_t;\psi_t),\\
    \epsilon \dot{P}_t &= -\gamma P_t + \nabla_x H(X_t,P_t,U_t,\beta_t;\psi_t) + D_x \psi_t(X_t)\dot{X}_t,\\
    U_t &= -G(X_t)^{\rm T}P_t
    \end{align}
\end{subequations}
since $1+\epsilon\approx 1$ for sufficiently small $\epsilon$. Here we have dropped the dependence of $X_t$, $P_t$, $U_t$, and $\beta_t$ on labels $a$.

Although it has been shown in \cite{R25b} that (\ref{eq:identity}) holds for any choice of $\beta_t$, it is often desirable that the law of $X_t$, as defined by (\ref{eq:FSDE}) agrees with the law of $X_t$ defined by (\ref{eq:PMP}a). This requirement leads to
\begin{equation} \label{eq:natural beta}
    \beta_t^{\rm SOC}  = \frac{1}{2}\Sigma \nabla_x \log \pi_t(X_t)
\end{equation}
with $\pi_t$ denoting the law of $X_t$ as defined by (\ref{eq:PMP}a).

The functions $v_\ast(x)$ and $\psi_\ast(x)$ are derived under the assumption of fully observable twins. This limitation can be resolved by a formal application of the separation or equivalence principle \cite{Handel07,WW81}; {\it i.e.},
\begin{equation} \label{eq:control equivalence}
U_t = -\int_{\mathbb{R}^{d_x}} G^{\rm T}(x) \psi_\ast(x) \,\pi_t(x|\{Y_{t_n\le t}^\dagger\}){\rm d}x
\end{equation}
provides an open loop control for partially observed diffusion processes.
However, it is well-known that (\ref{eq:control equivalence}) is suboptimal for partially observed nonlinear diffusion processes. Furthermore, it is desirable to maintain the online learning character of the evolution equations (\ref{eq:PMP}). In the next section, we develop an extension of (\ref{eq:constraint})--(\ref{eq:Hamiltonian density}) to partially observed physical twins. The proposed methodology  combines the ensemble Kalman filter formulations of Section \ref{sec:EnKF} with the McKean--Pontryagin formulation of this section.

%%%%%%%%%%%%%%%%%%%%%%%%%%%%%%%%%%%%%%%%%%%%%%%%%%%
%
\section{McKean--Pontryagin for Partially Observed Twins}
\label{sec:McKean--Pontryagin}
%
%%%%%%%%%%%%%%%%%%%%%%%%%%%%%%%%%%%%%%%%%%%%%%

In this section, we extend (\ref{eq:PMP}) to partially observed physical twins. We again consider both cases; time-continuous (\ref{eq:cobs}) and time-discrete observations (\ref{eq:obs}). We start with
the time-continuous case for ease of presentation.

%%%%%%%%%%%%%%%%%%%%%%%%%%%%%%%%%%%%%%%%%%%%%%%%%%%%%%%
\subsection{Continuously Observed Physical Twins}
%%%%%%%%%%%%%%%%%%%%%%%%%%%%%%%%%%%%%%%%%%%%%%%%%%%%%%

The first step is to replace the stochastic differential equation (\ref{eq:FSDE}) with the data driven formulation (\ref{eq:EnKF}). Furthermore, the desired control $U_t$ depends on the law $\pi^{\rm KBF}_t$ of $X_t$ as defined by (\ref{eq:EnKF}). In other words, instead of closed loop control (\ref{eq:optimal control law}), we will derive evolution equations for open loop control, which we denote by $U_t^\ast$. 

The open loop control $U_t^\ast$ is independent of the labels $a$. Hence, we consider the total Hamiltonian
\begin{equation} \label{eq:total Hamiltonian}
\mathcal{H}(X_t,P_t,u,\beta_t,\psi_t) = \int_{\mathbb{R}^{d_x}} H(X_t(a),P_t(a),u,\beta_t(a);\psi_t)\,\pi_0(a){\rm d}a
\end{equation}
with Hamiltonian density (\ref{eq:Hamiltonian density}) defined as before. Indeed, setting the gradient of (\ref{eq:total Hamiltonian}) with respect to $u$ to zero leads to the open loop control 
\begin{equation} \label{eq:openU}
    U_t^\ast = -\int_{\mathbb{R}^{d_x}} G(X_t(a))^{\rm T}P_t(a)\,\pi_0(a){\rm d}a .
\end{equation}

Next, we need to choose the drift term $\beta_t(a)$ in (\ref{eq:total Hamiltonian}). The definition of open loop control (\ref{eq:openU}) requires that the law $\pi_t$ of $X_t$ coincides with the law $\pi_t^{\rm KBF}$ of the mean-field Kalman--Bucy filter process (\ref{eq:EnKF}), which leads to
\begin{equation} \label{eq:beta_DT}
    \beta_t^{\rm DT}(a) = \beta_t^{\rm SOC}(a) - C_t^{xh} R^{-1} \left(\dot{Y}_t^\dagger - \frac{1}{2} \left(h(X_t(a)) + m^h_t \right)  \right).
\end{equation} 

Finally, after taking gradients of the Hamiltonian density (\ref{eq:Hamiltonian density}) with respect to $x$ and $p$, the McKean--Pontryagin formulation (\ref{eq:PMP}) is now replaced by the evolution equations 
\begin{subequations} \label{eq:PMPopen}
    \begin{align}
    \dot{X}_t(a) =& \,\,\nabla_p H(X_t(a),P_t(a),U_t^\ast,\beta^{\rm DT}_t(a);\psi_t) ,\\
    \epsilon \dot{P}_t(a) =& \,\,-\gamma P_t(a) + \nabla_x H(X_t(a),P_t(a),U_t^\ast,\beta_t^{\rm DT}(a);\psi_t) \nonumber \\
    & \,\, + (1+\epsilon) D_x \psi_t(X_t(a))\dot{X}_t(a),\\
    0=&\,\nabla_u \mathcal{H}(X_t,P_t,U_t^\ast,\beta_t,\psi_t)
    \end{align}
\end{subequations}
subject to (\ref{eq:constraint}). The initial conditions are $X_0(a) = a \sim \pi_0$ and $P_0(a) = 0$. Equations (\ref{eq:PMPopen}) are of mean-field type through $\psi_t (x)$, which is determined by the implicit relation (\ref{eq:constraint}), the function $\beta_t^{\rm DT}$, which involves the law $\pi_t=\pi_t^{\rm KBF}$ of $X_t$ and expectation values with respect to that law, and the control (\ref{eq:openU}), which now satisfies
\begin{equation} \label{eq:openU2}
    U_t^\ast = -\int_{\mathbb{R}^{d_x}} G(x)^{\rm T}\psi_t(x)\,\pi_t(x){\rm d}x .
\end{equation}

Upon ignoring labels $a \in \mathbb{R}^{d_x}$ from now on, the McKean--Pontryagin evolution equations (\ref{eq:PMPopen}a)-(\ref{eq:PMPopen}b) become
\begin{subequations} \label{eq:Hamiltonian ODE continuous}
    \begin{align}
        \dot{X}_t =&\,\, b(X_t) + G(X_t)U_t^\ast - \frac{1}{2}\Sigma \nabla_x \log \pi_t(X_t)  \nonumber \\
        &\,\,+C_t^{xh} R^{-1} \left(\dot{Y}_t^\dagger - \frac{1}{2} \left(
h(X_t) + m^h_t \right) \right) ,\\
        \epsilon \dot{P}_t 
        =&\,\, -\gamma P_t + \left(D_x \left(b(X_t)+G(X_t)U_t^\ast\right) \right)^{\rm T}P_t \nonumber \\ & \,\,+ \nabla_x c(X_t) + \mathcal{L}_{\beta^{\rm SOC}_t}  \psi_t(X_t)\nonumber \\
        &\,\, 
        +  D_x \psi_t(X_t)\left((1+\epsilon)\dot{X}_t + \beta_t^{DT}(X_t)-\beta_t^{\rm SOC}(X_t)\right)  .
    \end{align}
\end{subequations}
Here we have exploited (\ref{eq:gradient hjb}), which implies that the Jacobian $D_x \psi_t(x)$ is symmetric. Furthermore, we have introduced the generator 
$\mathcal{L}_{\beta_t^{\rm SOC}}$ via
\begin{equation} \label{eq:generator}
\mathcal{L}_{\beta_t^{\rm SOC}} g = (\beta_t^{\rm SOC})^{\rm T} \nabla_x g 
+ \frac{1}{2} \Sigma : D_x^2 g
\end{equation}
for scalar-valued functions $g(x)$ with obvious component-wise generalization to vector-valued $\psi_t(x)$. We also note that $\beta_t^{\rm SOC} = \mathcal{L}_{\beta_t^{\rm SOC}}X_t$. See \cite{R25b} for more details.

Formulation (\ref{eq:Hamiltonian ODE continuous}) together with the control (\ref{eq:openU}) provides the desired digital twin for the partially observed physical twin (\ref{eq:FSDE}) with continuous time observations provided by (\ref{eq:cobs}).

Let us discuss the role of the time-scale parameter $\epsilon >0$ in (\ref{eq:Hamiltonian ODE continuous}). Taking the limit $\epsilon \to 0$, we obtain the following stationary equation in $\psi_t$, which we denote by $\hat \psi_t$ for clarity of presentation:
\begin{equation} \label{eq:stationary}
\gamma \hat \psi_t = \nabla_x \left( \left(b+GU_t^\ast \right)^{\rm T}\hat \psi_t  + c 
+ \frac{1}{2}\nabla_x \cdot(\Sigma  \hat \psi_t) \right).
\end{equation}
Here we have used $P_t = \hat \psi_t(X_t)$, replaced $X_t$ by $x$, and subsequently dropped the argument. Upon applying the stationarity assumption (\ref{eq:stationary}), the open loop (\ref{eq:openU2}) becomes
\begin{equation} \label{eq:starionary control}
    U_t^\ast = -\int G(x)^{\rm T} \hat \psi_t(x)\,\pi_t(x){\rm d}x 
\end{equation}
with $\pi_t$ the law of $X_t$. Note that (\ref{eq:stationary}) corresponds to the stationary
Hamilton--Jacobi--Bellman equation 
\begin{equation} \label{eq:HJBs}
    \gamma \hat v_t = \left(b+GU_t^\ast \right)^{\rm T}\nabla_x \hat v_t
    + \frac{1}{2}\Sigma : D_x^2 \hat v_t + \frac{1}{2}\|U_t^\ast\|^2 + c.
\end{equation}
for the value function $\hat v_t(x)$ and $\hat \psi_t = \nabla_x \hat v_t$; compare (\ref{eq:HJB}). However, the optimal control is now determined by (\ref{eq:starionary control}) instead of the closed loop control (\ref{eq:control equivalence}), which is based on the stationary solution $v_\ast$ of (\ref{eq:HJB}) for fully observed diffusion processes and $\psi_\ast = \nabla_x v_\ast$  \cite{R25b}. Note that a change in $\pi_t$ implies a modification in $\hat \psi_t$ and the combination of both leads to open loop control (\ref{eq:starionary control}). 

The quasi-stationary value function $\hat v_t$ takes full account of the current filtering distribution $\pi_t$ while ignoring the impact of future observations. Hence, our formulation still only provides  an approximation to the stationary measure-valued value function, which arises from partially observed Markov decision processes \cite{Bensoussan92,Nisio15}. See Remark \ref{rem1} below. However, we note that our formulation becomes exact in the limit $\epsilon \to 0$ if the separation principle applies \cite{Handel07,WW81}; {\it i.e.},
\begin{equation}
    U_t^\ast = -\int G(x)^{\rm T} \hat \psi_t(x)\,\pi_t(x){\rm d}x = 
    -\int G(x)^{\rm T} \psi_\ast (x)\,\pi_t(x){\rm d}x.
\end{equation}

In practice, one would choose $\epsilon$ not too small in (\ref{eq:Hamiltonian ODE continuous}) since small values of $\epsilon > 0$ lead to stiff differential equations; but also not too large for (\ref{eq:Hamiltonian ODE continuous}) to stay close to the quasi-equilibrium solution $\hat \psi_t$ determined by (\ref{eq:stationary}). Furthermore, (\ref{eq:PMPopen}) can be simplified for $\epsilon$ sufficiently small as already done in the fully observed case; {\it i.e.}~(\ref{eq:PMPs}).

\begin{remark} \label{rem1}
The Hamiltonian density (\ref{eq:Hamiltonian density}) does not yet take into account the impact of future observations. A possible modification is provided by the Hamiltonian density
%\begin{subequations}
\begin{align}
    \tilde H(x,p,u,\beta;\psi,\pi) =& \,H(x,p,u,\beta;\psi) -\frac{1}{2}p^{\rm T}C^{xh}R^{-1}(h(x)-m^h) \nonumber \\
    &\,+\, \frac{1}{2} (C^{xh}R^{-1}C^{hx}): D_x\psi(x).
\end{align}
%\end{subequations}
Since the law of $X_t$ should still agree with the law $\pi_t^{\rm KBF}$, $\beta_t$
needs to be adjusted accordingly; {\it i.e.},
    \begin{equation} \label{eq:beta_DT2}
    \tilde \beta_t^{\rm DT}(a) = \beta_t^{\rm SOC}(a) - C_t^{xh} R^{-1} \left(\dot{Y}_t^\dagger - m^h_t   \right).
\end{equation}
The associated stationary Hamilton--Jacobi--Bellman equation (\ref{eq:HJBs}) is given by
%\begin{subequations}
    \begin{align}
    \gamma \hat v_t =&\,\, \left(b+GU_t^\ast - \frac{1}{2}C_t^{xh}R^{-1}\left(
    h-m_t^h\right)\right)^{\rm T}\nabla_x \hat v_t \nonumber \\
    &\,\,
    + \frac{1}{2}\left(\Sigma+C_t^{xh}R^{-1}C_t^{hx}\right) : D_x^2 \hat v_t + \frac{1}{2}\|U_t^\ast\|^2 + c.
\end{align}
%\end{subequations}
A detailed exploration of this and other formulations relative to optimal control \cite{PM13} implied by the underlying partially observed Markov decision process is left for future research.
\end{remark}

%%%%%%%%%%%%%%%%%%%%%%%%%%%%%%%%%%%%%%%%%%%%%%%%%%%
\subsection{Discretely Observed Physical Twins}
%%%%%%%%%%%%%%%%%%%%%%%%%%%%%%%%%%%%%%%%%%%%%%%%%%%

The general structure of the evolution equations in the states $X_t$ and co-states $P_t$ remains largely unchanged under discrete-time observations (\ref{eq:obs}). In particular, open loop control (\ref{eq:openU}) and the Hamiltonian density (\ref{eq:Hamiltonian density}) remain as before. The only change arises from the modified data driven term in the ensemble Kalman formulation (\ref{eq:mollied EnKF}), which results in 
\begin{equation}
\beta_t^{\rm DT}(X_t) = \beta_t^{\rm SOC}(X_t) - \sum_{n\ge 1} C_t^{xh} R^{-1} \left(\dot{Y}_{t_n}^\dagger - \frac{1}{2} \left(
h(X_t) + m^h_t \right) \right)\,\delta_{t_n}(t).
\end{equation}
The resulting McKean--Pontryagin evolution equations become
\begin{subequations} \label{eq:Hamiltonian ODE discrete}
    \begin{align}
        \dot{X}_t =&\,\, b(X_t) + G(X_t)U_t^\ast - \frac{1}{2}\Sigma \nabla_x \log \pi_t(X_t)  \nonumber \\
        &\,\,+\sum_{n\ge 1} C_t^{xh} R^{-1} \left(\dot{Y}_{t_n}^\dagger - \frac{1}{2} \left(
h(X_t) + m^h_t \right) \right) \delta_{t_n}(t),\\
        \epsilon \dot{P}_t 
        =&\,\, -\gamma P_t + \left(D_x \left(b(X_t)+G(X_t)U_t^\ast\right) \right)^{\rm T}P_t  \nonumber \\
        &\,\, + \nabla_x c(X_t) + \mathcal{L}_{\beta_t^{\rm SOC}}  \psi_t(X_t) \nonumber\\
        & \,\,+  D_x \psi_t(X_t)\left((1+\epsilon) \dot{X}_t + \beta_t^{DT}(X_t)-\beta_t^{\rm SOC}(X_t)\right)
    \end{align}
\end{subequations}
and the law of $X_t$ is equal to $\pi_t^{\rm KF}$ for given control $U_s^\ast$, $s\le t$.

Formulation (\ref{eq:Hamiltonian ODE discrete}) together with the control (\ref{eq:openU}) provides the desired digital twin for the partially observed physical twin (\ref{eq:FSDE}) with discrete-time observations provided by (\ref{eq:obs}).

%%%%%%%%%%%%%%%%%%%%%%%%%%%%%%%%%%%%%%%%%%%
%
\section{Algorithmic implementation} \label{sec:implementation}
%
%%%%%%%%%%%%%%%%%%%%%%%%%%%%%%%%%%%%%%%%%

We introduce $M$ particles with states $X_t^{(i)} \in \mathbb{R}^{d_x}$ and co-states $P_t^{(i)} \in \mathbb{R}^{d_x}$, $i=1,\ldots,M$. Based on the Schr\"odinger bridge approach taken in \cite{R25b}, we obtain the following interacting particle approximation of (\ref{eq:Hamiltonian ODE continuous}):
\begin{subequations} \label{eq:interacting particle formulation}
    \begin{align}
        \dot{X}_t^{(i)} =&\,\, b(X_t^{(i)}) + G(X_t^{(i)})U_t^\ast
        - \sum_{j=1}^M \mu_t^{(ij)}X_t^{(j)} \nonumber \\
        & \,\,+C_t^{xh}R^{-1}\left( \dot{Y}_t^\dagger - \frac{1}{2}\left(
        h(X_t^{(i)}) + m_t^h \right) \right),\\
        \epsilon \dot{P}_t^{(i)} =&\,\, -\gamma P_t^{(i)}+ \left( D_x \left( b(X_t^{(i)})+ G(X_t^{(i)}U_t^\ast \right)\right)^{\rm T}P_t^{(i)} +\nabla_x c(X_t^{(i)}) \nonumber
        \\
        & \,\,+ \sum_{j=1}^M \mu_t^{(ij)}P_t^{(j)} + (1+\epsilon) 
        D_x \psi_t(X_t^{(i)}) \dot{X}_t^{(i)} \nonumber \\
        & \,\,-D_x \psi_t(X_t^{(i)}) C_t^{xh}R^{-1}\left( \dot{Y}_t^\dagger - \frac{1}{2}\left(
        h(X_t^{(i)}) + m_t^h \right) \right)
    \end{align}
\end{subequations}
with control
\begin{equation} \label{eq:empirical control}
    U_t^\ast = \frac{1}{M}\sum_{i=1}^M G(X_t^{(i)})^{\rm T}P_t^{(i)},
\end{equation}
empirical mean approximations
\begin{equation}
m_t^x = \frac{1}{M}\sum_{i=1}^M X_t^{(i)}, \qquad 
m_t^h = \frac{1}{M}\sum_{i=1}^M h(X_t^{(i)}),
\end{equation}
and empirical covariance matrix
\begin{equation}
    C_t^{xh} = \frac{1}{M}\sum_{i=1}^M (X_t-m_t^x)(h(X_t)-m^h_t)^{\rm T}.
\end{equation}
The function $\psi_t(x)$ is approximated using Nadaraya--Watson kernel regression \cite{Bierens}; {\it i.e.},
\begin{equation} \label{eq:NWKR}
\psi_t(x):=
\frac{\sum_{i=1}^M \exp \left(-\frac{1}{2\delta}\|x - X_t^{(i)}\|^2 \right) P_t^{(i)}}{\sum_{i=1}^M \exp \left(-\frac{1}{2\delta}\|x - X_t^{(i)}\|^2 \right)}
\end{equation}
for a suitable parameter choice $\delta >0$, and
\begin{equation}
    D_x \psi_t (X_t) \,d \approx \frac{1}{\Delta t} \left(\psi_t(X_t+\Delta t\, d) - \psi_t(X_t)\right)
\end{equation}
for any vector $d \in \mathbb{R}^{d_x}$ and suitable step-size $\Delta t >0$. In order to approximate the drift terms arising from the generator (\ref{eq:generator}), we employ the Schr\"odinger bridge approximation of \cite{R25b}. Specifically, the coefficients $\mu_t^{(ij)} \in \mathbb{R}$, $i,j=1,\ldots,M$, in (\ref{eq:interacting particle formulation}) are defined by
\begin{equation}
    \mu_t^{(ij)} = \frac{1}{\alpha} \left( v_t^{(i)} d_t^{(ij)} v_t^{(j)} - \delta_{ij}\right), \qquad d_t^{(ij)} = 
    e^{-\frac{1}{2\alpha}\|X_t^{(i)}-X_t^{(j)}\|_\Sigma^2},
\end{equation}
with $v_t^{(i)}\ge 0$, $i=1,\ldots,M$, chosen such that 
\begin{equation}
\sum_i \mu_t^{(ij)} = \sum_j \mu_t^{(ij)} = 0.
\end{equation}
Here, $\delta_{ij}$ denotes the Kronecker delta function; {\it i.e.}, $\delta_{ij} = 1$ if $i=j$ and $\delta_{ij} = 0$ otherwise, and
\begin{equation}
    \|x\|_\Sigma^2 = x^{\rm T}\Sigma^{-1}x.
\end{equation}
The regularization parameter $\alpha>0$ should be chosen appropriately \cite{WR20}. A particularly elegant formulation arises from the choice $\alpha = \Delta t$ with $\Delta t>0$ the step-size of a Euler discretization of the equations of motion (\ref{eq:interacting particle formulation}) \cite{GLRY24}. In \cite{MarshallCoifman,WR20}, a fast iterative algorithm has been given to compute the coefficients $\{v_t^{(i)}\}$. 

The ensemble Kalman--Bucy filter is often implemented for ensemble sizes $M\ll d_x$ \cite{Evensenetal2022} leading to a computational complexity that scales linearly in $d_x$ and quadratically in $M$. To achieve computational robustness, such implementations require localization and inflation \cite{Evensenetal2022}. 
The same computational complexity arises from the proposed implementation of the McKean--Pontryagin control formulation using linear approximations
\begin{equation} \label{eq:linear approximation}
\nabla_x \log \pi_t(x) \approx -(C_t^{xx})^{-1}(x-m_t^x), \qquad
\Psi_t(x) \approx C_t^{xp}(x-m_t^x) + m_t^p.
\end{equation}
Here, $C_t^{xp}$ denotes the covariance matrix between $X_t$ and $P_t$ and $m_t^p$ the mean of $P_t$. See \cite{R25b} for a variational derivation of the approximations (\ref{eq:linear approximation}). 

Covariance localization \cite{Evensenetal2022,reich2015probabilistic} can be applied to (\ref{eq:linear approximation}) if $M\ll d_x$. We note that a localization strategy for Schr\"odinger bridges has also been introduced in \cite{gottwald2024localized} in the context of diffusion generative modeling. This localization strategy can be extended to McKean--Pontryagin control; {\it i.e.}, Nadaraya--Watson kernel regression and the Schr\"odinger bridge approximation in (\ref{eq:interacting particle formulation}).

The interacting particle system (\ref{eq:interacting particle formulation}) is discretized in time by the standard forward Euler method with step-size $\Delta t>0$; thus providing a fully discretized digital twin. We also use $\tau = \Delta t$ in (\ref{eq:obsc2}).

The same numerical approximations can be applied to formulation (\ref{eq:Hamiltonian ODE discrete}) for discrete-time observations (\ref{eq:obs}).

%%%%%%%%%%%%%%%%%%%%%%%%%%%%%%%%%%%%%%%
%
\section{Numerical examples} \label{sec:example}
%
%%%%%%%%%%%%%%%%%%%%%%%%%%%%%%%%%%%%%%

We present results from three numerical experiments. We start with a controlled Lorenz-63 system as proposed in \cite{KK24}. The second example considers the classical inverted pendulum control problem \cite{Meyn}. The third example is based on the chaotic $d_x =40$ dimensional Lorenz-96 model \cite{lorenz96}, which is partially observed and controlled with the aim of stabilizing $x=0$. All three cases lead to physical twins in the form of ordinary differential equations; {\it i.e.}, $\Sigma = 0$ in (\ref{eq:FSDE}). The time-scale parameter $\epsilon$ in (\ref{eq:interacting particle formulation}) is set to one for the first two numerical experiments and to $\epsilon = 0.1$ for the rapidly changing and highly chaotic Lorenz-96 system. 

%%%%%%%%%%%%%%%%%%%%%%%%%%%%
%
\subsection{Controlled Lorenz-63 system}
%
%%%%%%%%%%%%%%%%%%%%%%%%%%%%%

\begin{figure}
\begin{center}
\includegraphics[width=0.47\textwidth]{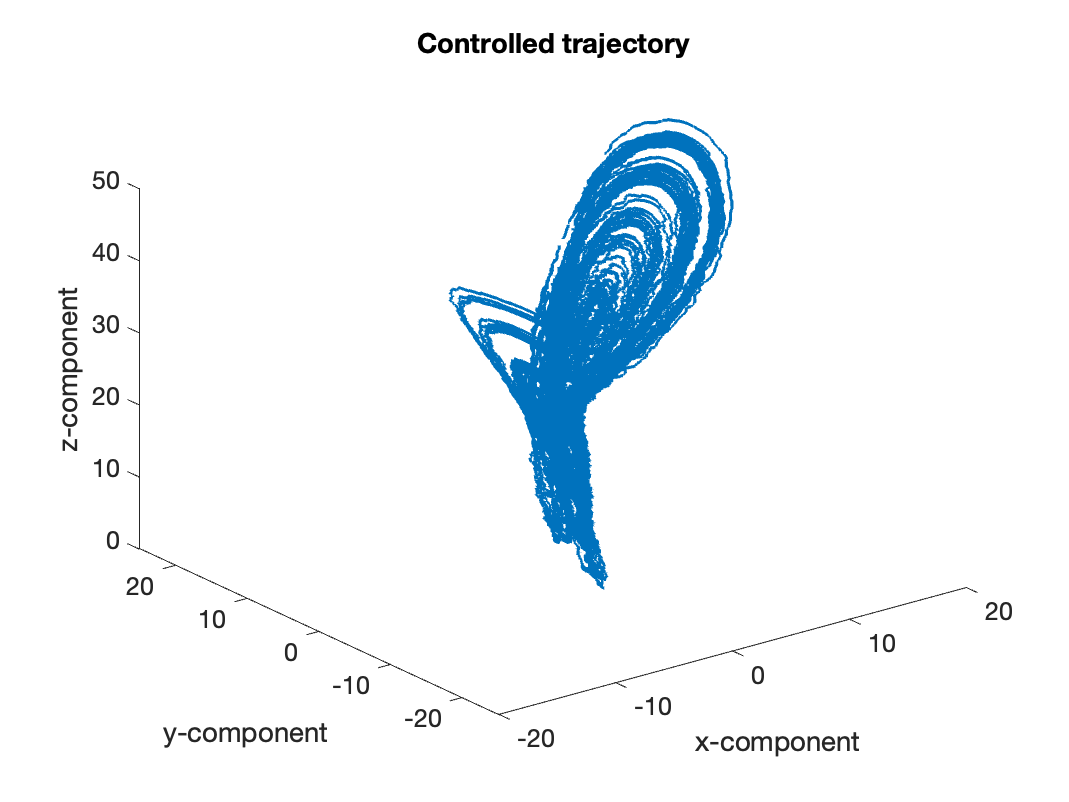}
$\quad$ 
\includegraphics[width=0.47\textwidth]{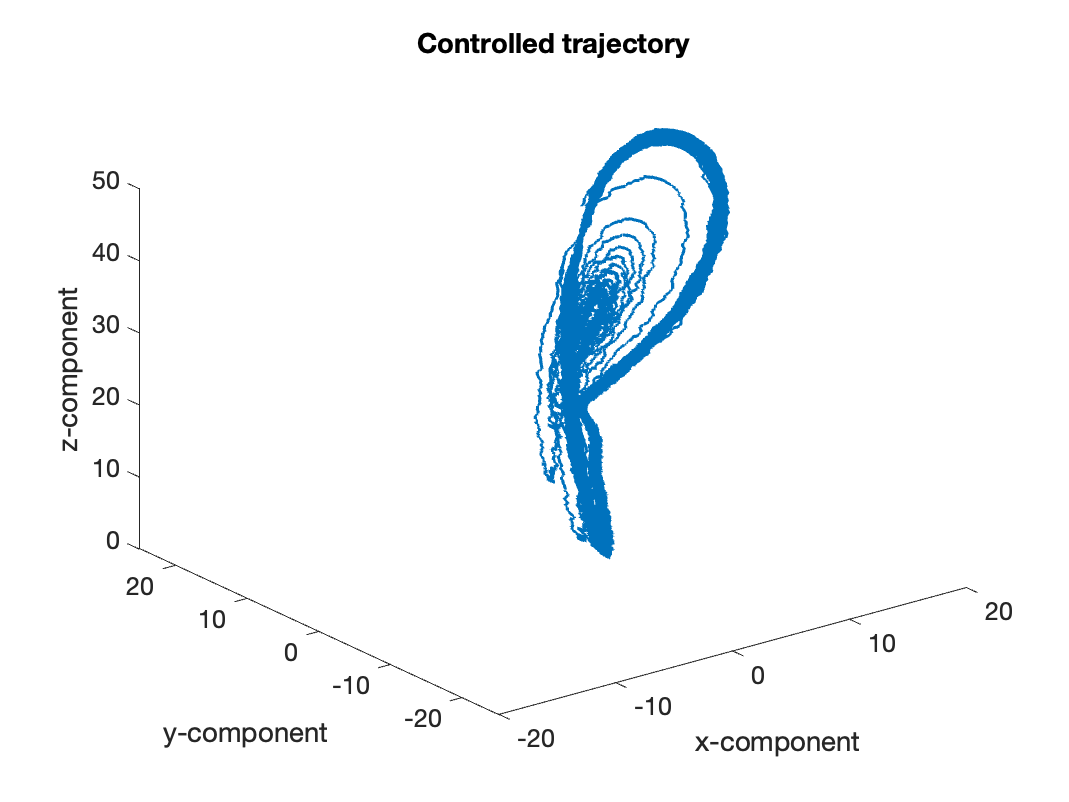}
\end{center}
\caption{Controlled Lorenz-63 model.~Displayed is the
three-dimensional trajectory of the  mean $\{m^x_t\}$ over the time interval $t \in [0,100]$.  Left panel: control restricted to range $|U_t|\le 50$; right panel: control restricted to range $|U_t|\le 100$. After an initial transient, the trajectory enters a quasi-period orbit for the larger threshold value. The smaller threshold value still allows for some transitions to negative ${\rm x}_t$ values and the chaotic nature of the Lorenz-63 system is retained.} 
\label{fig1}
\end{figure}

\begin{figure}
\begin{center}
\includegraphics[width=0.47\textwidth]{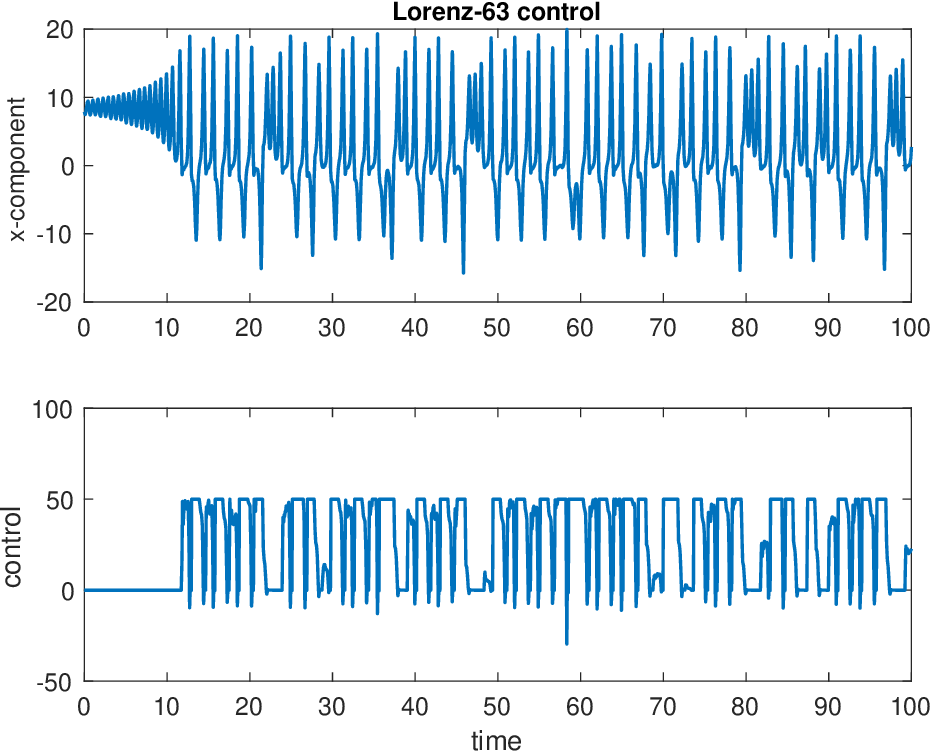}
$\quad$
\includegraphics[width=0.47\textwidth]{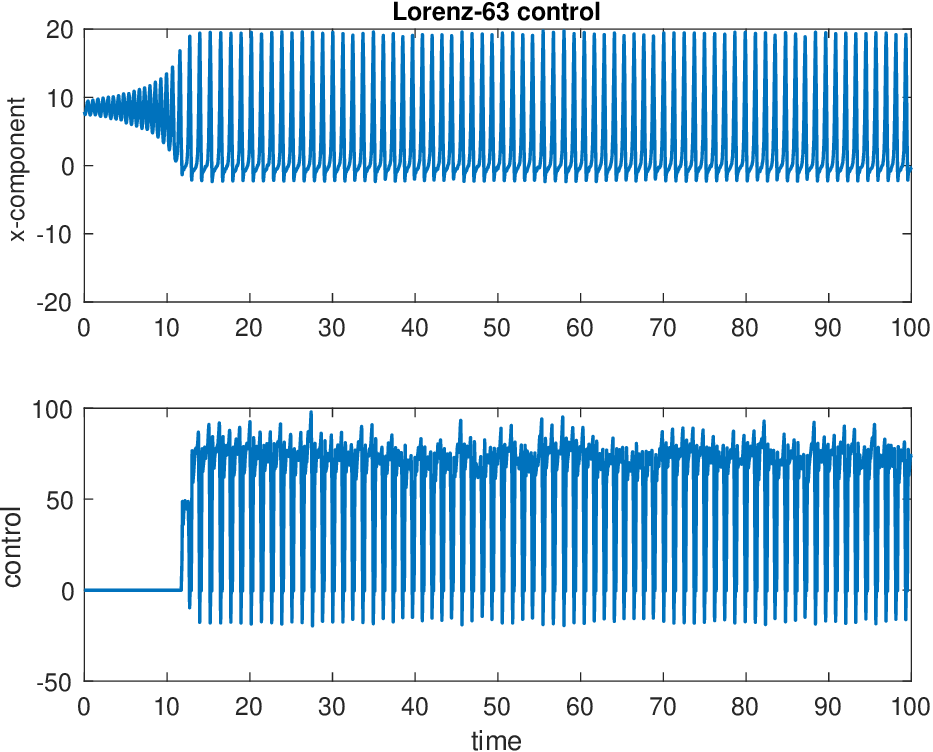}
\end{center}
\caption{Controlled Lorenz-63 model.~First component of the mean $m_t^x \in \mathbb{R}^3$ and associated control term $U_t$ as a function of time. Left panels: control restricted to range $|U_t|\le 50$; right panels: control restricted to range $|U_t|\le 100$.} 
\label{fig2}
\end{figure}

\begin{figure}
\begin{center}
\includegraphics[width=0.47\textwidth]{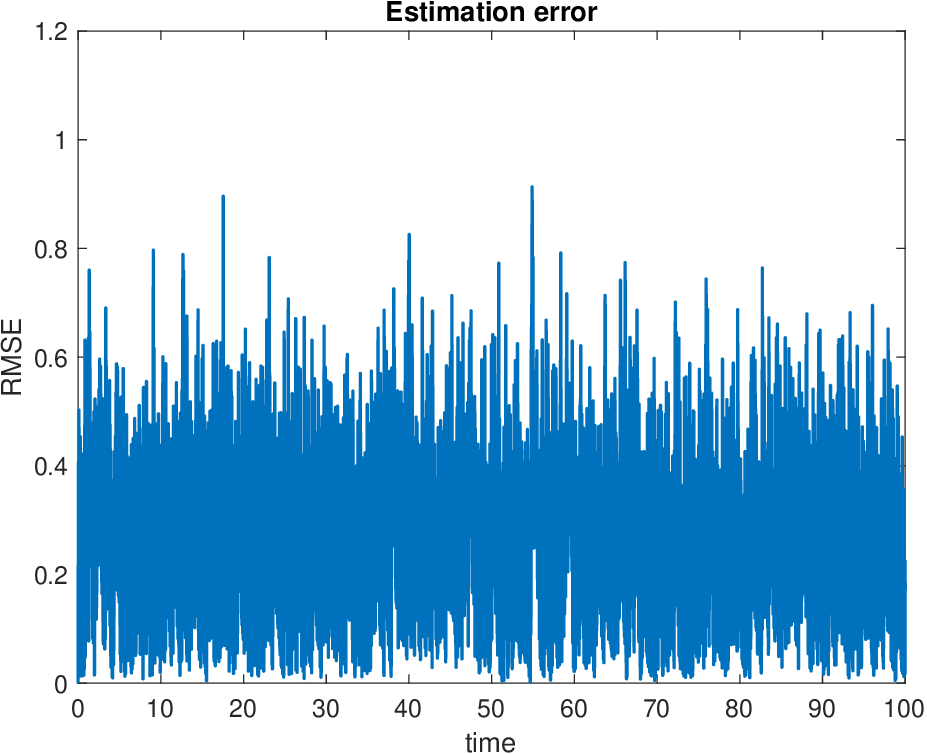}
$\quad$
\includegraphics[width=0.47\textwidth]{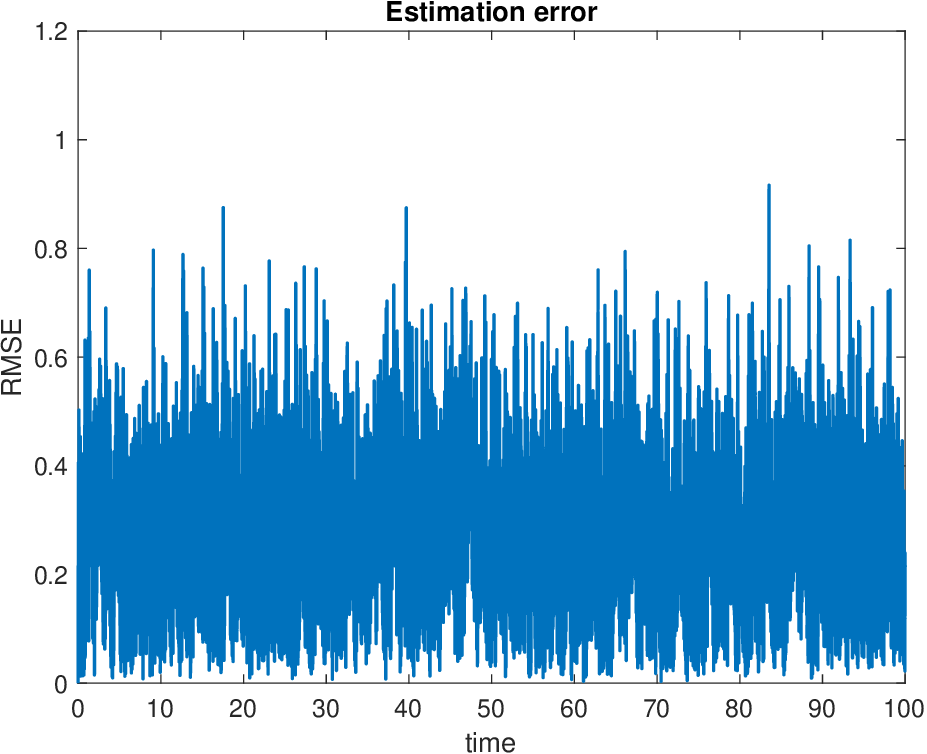}
\end{center}
\caption{Controlled Lorenz-63 model.~ Estimation error; {\it i.e.}~root mean square error between physical twin states $X_t^\dagger$ and digital twin mean states $m_t^x$, as a function of time. Left panel: control restricted to range $|U_t|\le 50$; right panel: control restricted to range $|U_t|\le 100$.} 
\label{fig3}
\end{figure}

We follow \cite{KK24} and consider an optimal control problem for the Lorenz-63 model
\begin{equation} \label{eq:l63}
\dot{X}_t = b(X_t) + G U_t, \qquad b(x) = \left( \begin{array}{l} \sigma({\rm y}-{\rm x})\\
-{\rm x}{\rm z} + r{\rm x}-{\rm y}\\ {\rm x}{\rm y}-b{\rm z} \end{array}\right),
\end{equation}
in the state variable $x = ({\rm x},{\rm y},{\rm z})^{\rm T} \in \mathbb{R}^3$ and with
parameters $\sigma =10$, $r=28$, and $b = 8/3$ \cite{lorenz1963deterministic}. The scalar-valued control $U_t$ acts only on the ${\rm x}$-component; {\it i.e.}~$G = (1,0,0)^{\rm T}$. The control problem is to restrict the solutions of (\ref{eq:l63}) to ${\rm x}_t\ge 0$ and the authors of \cite{KK24} introduce the running cost 
\begin{equation}
c(x) = \frac{\rho}{2} \left( \min({\rm x},0) \right)^2.
\end{equation}
In our experiments, we have followed this cost function and have set $\rho = 5000$. Contrary to the setting of \cite{KK24}, we have implemented a discounted cost function (\ref{eq:cost}) with $\gamma = 10$ instead of a receding horizon model predictive control approach \cite{Allgower,MPC}. 

Time-continuous observations (\ref{eq:cobs}) are obtained by solving the Lorenz system (\ref{eq:l63}); {\it i.e.}~the physical twin, with initial condition 
\begin{equation}
X_0^\dagger = (7.8590,7.1136,27.2293)^{\rm T}
\end{equation}
and control $U_t^\ast$ provided by its digital twin. We observe the second component
of $X_t^\dagger \in \mathbb{R}^3$ subject to Brownian noise with variance 
$R = 0.0025$; {\it i.e.}~$h(x) = {\rm y} \in \mathbb{R}$. Observing the first component only leads to similar results; but for a smaller measurement error variance $R = 0.001$.

Since the forward operator is linear; {\it i.e.} $h(x) = Hx$, we use the following time-dependent approximation to the covariance matrix $C_t^{xh} = C_t^{xx}H^{\rm T}$
with
\begin{equation} \label{eq:inflation}
C_t^{xx} = \frac{1}{M} \sum_{i=1}^M (X_t^{(i)}-m_t^x)(X_t^{(i)}-m_t^x)^{\rm T}
+ \sigma I
\end{equation}
The inflation parameter is set to $\sigma = 0.2$ This approximation corresponds to additive ensemble inflation widely used in ensemble Kalman filtering \cite{law2015data,reich2015probabilistic,Evensenetal2022}.

In terms of the digital twin for (\ref{eq:l63}), we have implemented the Schr\"odinger bridge based McKean--Pontryagin formulation (\ref{eq:interacting particle formulation}) with $\alpha = 0.1$, step-size $\Delta t = 0.001$, and $M=100$ and $M=4$ particles, respectively. The parameter $\delta$ in (\ref{eq:NWKR}) is set to $\delta = 1$. The initial co-states are $P_0^{(i)} = 0$. The numerical parameters are identical to those used in \cite{R25b} for the fully observed Lorenz-63 optimal control problem.

We added artificial diffusion (viscosity) via $\Sigma = 0.5I$ in (\ref{eq:Hamiltonian density}) to regularize the underlying Hamilton--Jacobi--Bellman equation (\ref{eq:HJBs}). This regularization corresponds to the viscosity solution approach of optimal control \cite{CL83,Carmona}.

The results with $M=100$ particles can be found in Figures \ref{fig1}, \ref{fig2}, and \ref{fig3}. The impact of the control can be clearly detected in Figures \ref{fig1} and \ref{fig2}, where we present the results of two different scenarios by enforcing that the control (\ref{eq:empirical control}) does not exceed a threshold value $U_{\rm max}$ in absolute value. The panels on the left display the results for $U_{\rm max} = 50$ while the panels on the right are for the larger threshold of $U_{\rm max} = 100$. As can be clearly seen, the smaller threshold value leads to some transitions to negative ${\rm x}_t$ values. The higher threshold value, on the other hand, nearly perfectly enforces ${\rm x}_t \ge 0$ thus eliminating the chaotic nature of the Lorenz-63 model. The initial spiral towards the quasi-periodic attractor is due to the chosen initialization 
\begin{equation}
X_0^{(i)} = X_0^\dagger + \Xi_0^{(i)}, \qquad \Xi_0^{(i)} \sim {\rm N}(0,0.1I),
\end{equation}
of the Lorenz-63 problem. Our results for $U_{\rm max} = 100$ are qualitatively in agreement with the findings of \cite{KK24},  which are based on a combination of data assimilation \cite{reich2015probabilistic,Evensenetal2022} and model predictive control \cite{Allgower,MPC}, as well as \cite{R25b}, which are based on the McKean--Pontryagin formulation for fully observed physical twins. Figure \ref{fig3} displays the root mean square error 
\begin{equation}
    \mbox{RMSE}_t := \sqrt{ \frac{1}{3} \|X_t^\dagger - m_t^x\|}
\end{equation}
between the trajectory of the physical twin, $X_t^\dagger$, and the mean state of its digital twin, $m_t^x$, which remains bounded throughout the simulation interval regardless of the value of $U_{\rm max}$. 

\begin{figure}
\begin{center}
\includegraphics[width=0.47\textwidth]{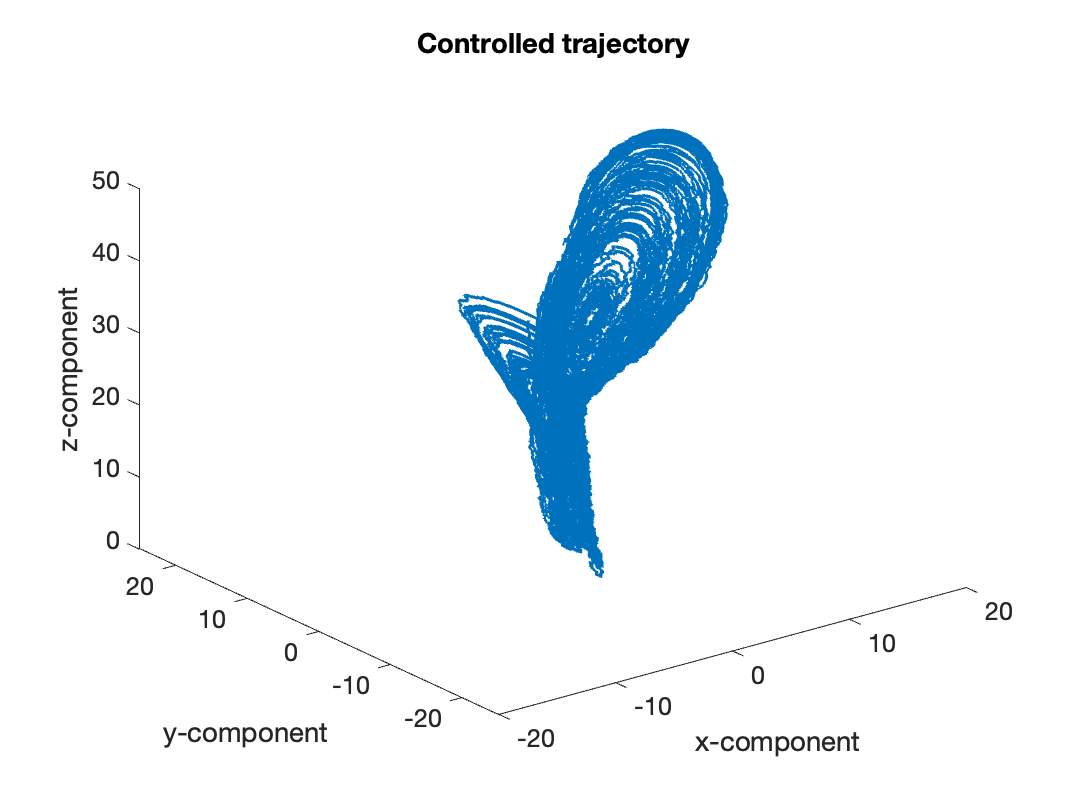}
$\quad$ 
\includegraphics[width=0.47\textwidth]{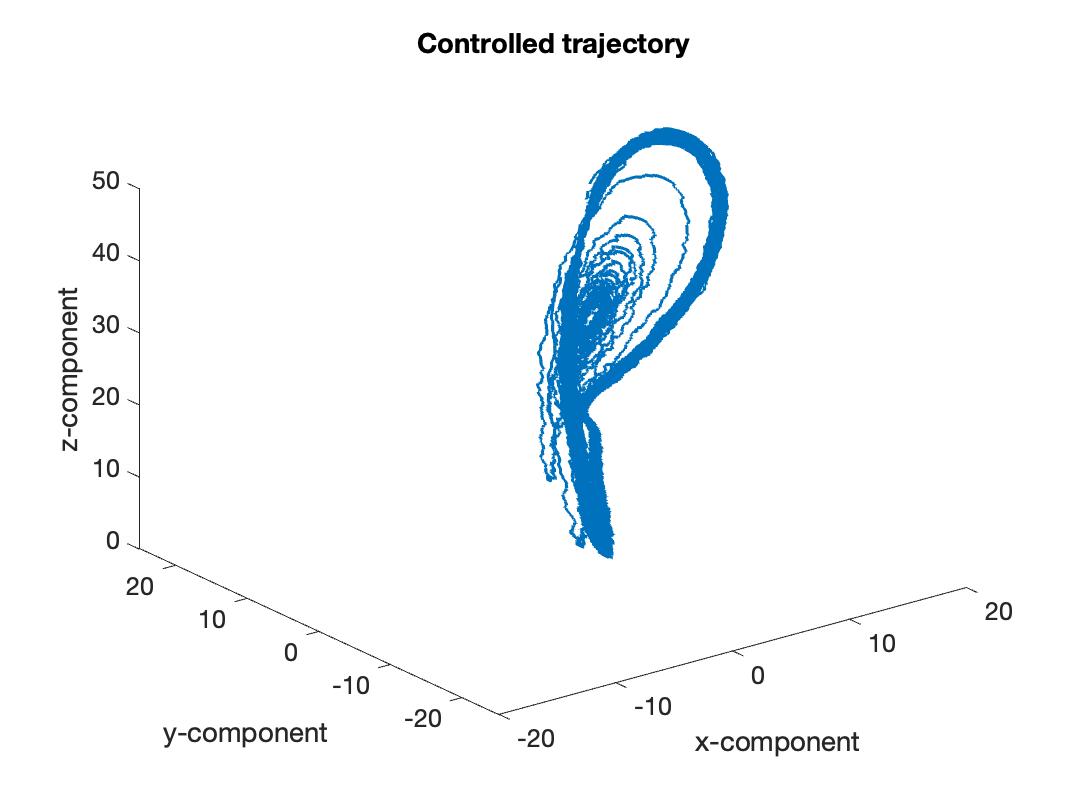}
\end{center}
\caption{Controlled Lorenz-63 model.~Same as Figure \ref{fig1} but for ensemble size 
$M=4$.} 
\label{fig4}
\end{figure}

\begin{figure}
\begin{center}
\includegraphics[width=0.47\textwidth]{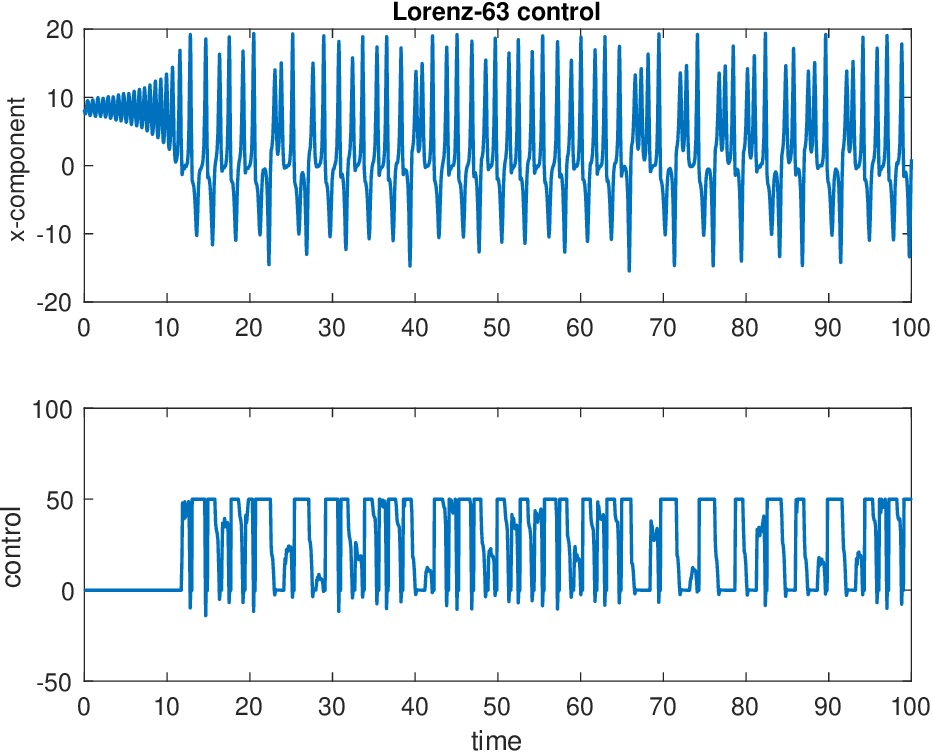}
$\quad$
\includegraphics[width=0.47\textwidth]{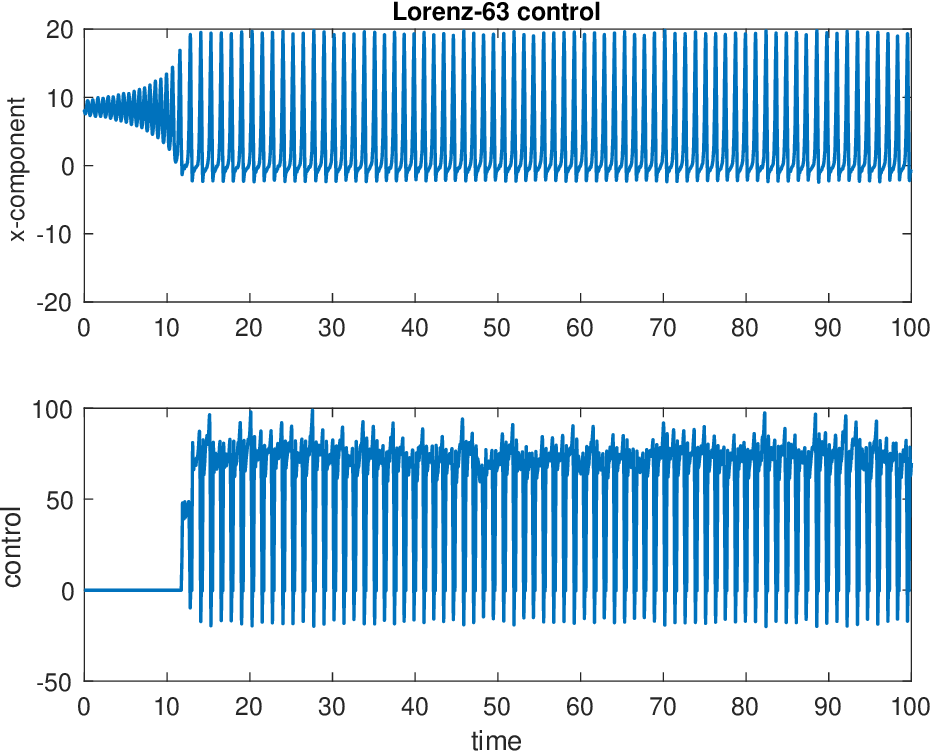}
\end{center}
\caption{Controlled Lorenz-63 model.~Same as Figure \ref{fig2} but for ensemble size 
$M=4$.} 
\label{fig5}
\end{figure}

\begin{figure}
\begin{center}
\includegraphics[width=0.47\textwidth]{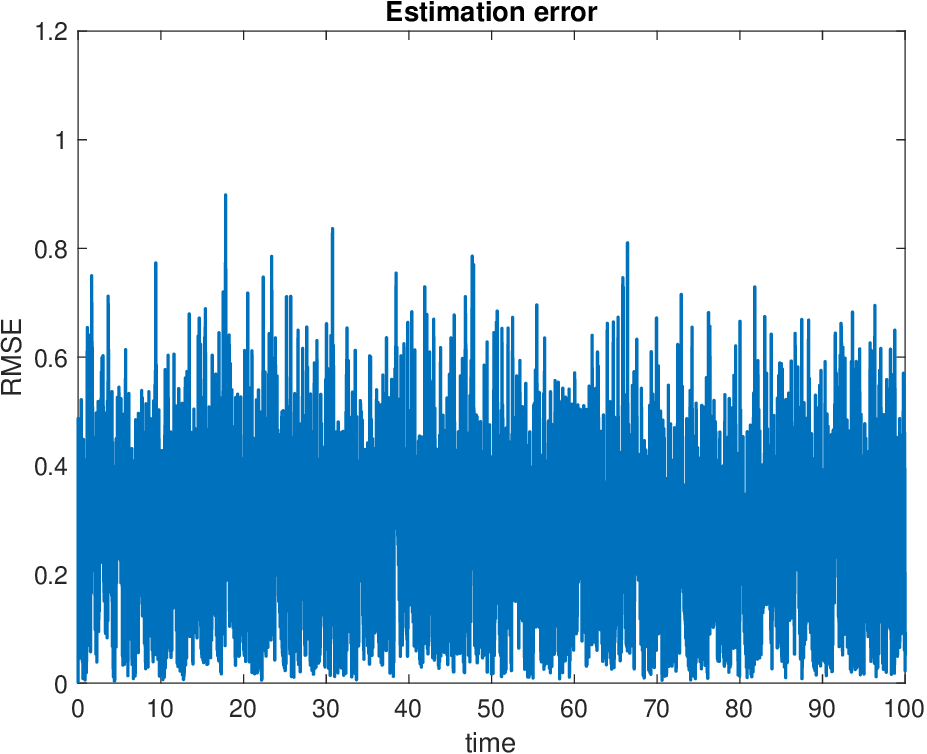}
$\quad$
\includegraphics[width=0.47\textwidth]{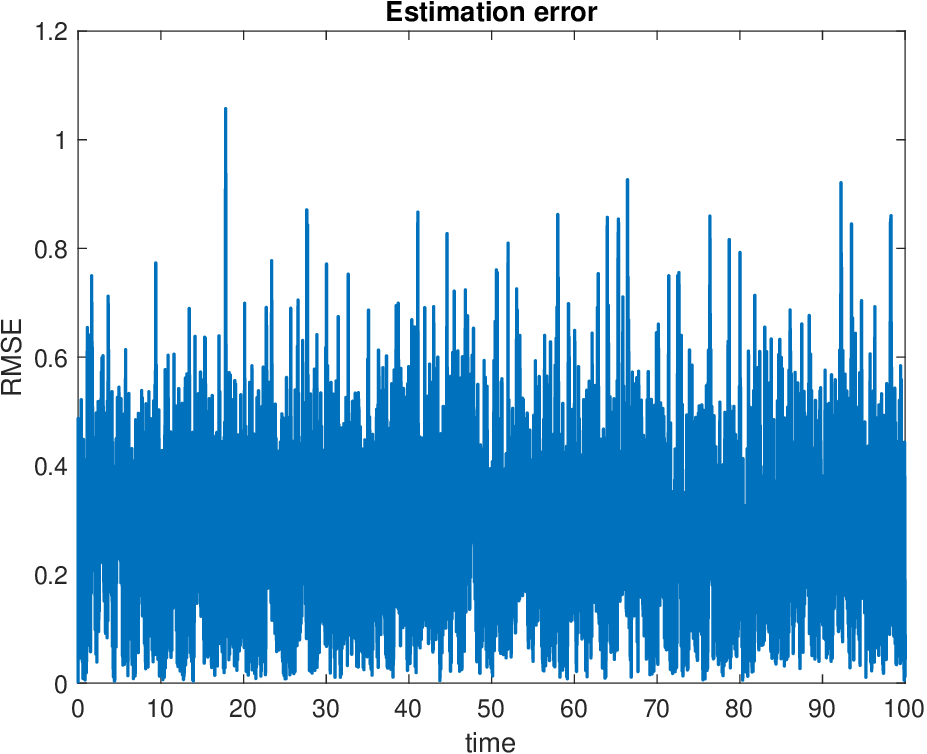}
\end{center}
\caption{Controlled Lorenz-63 model.~Same as Figure \ref{fig3} but for ensemble size 
$M=4$.} 
\label{fig6}
\end{figure}

In order to investigate the computational robustness of our approach in terms of ensemble sizes, in Figures \ref{fig4}, \ref{fig5}, and \ref{fig6} we show the corresponding results for $M=4$. We remind the reader that additive inflation is used in (\ref{eq:inflation}) which avoids $C_t^{xx}$ from becoming singular. The results agree well with those for $M=100$ and demonstrate the desirable computational robustness of the proposed methodology.

%%%%%%%%%%%%%%%%%%%%%%%%%%%%%%%%%
%
\subsection{Inverted pendulum} \label{sec:num_IP}
%
%%%%%%%%%%%%%%%%%%%%%%%%%%%%%%%%

\begin{figure}
\centerline{\includegraphics[width=0.47\textwidth]{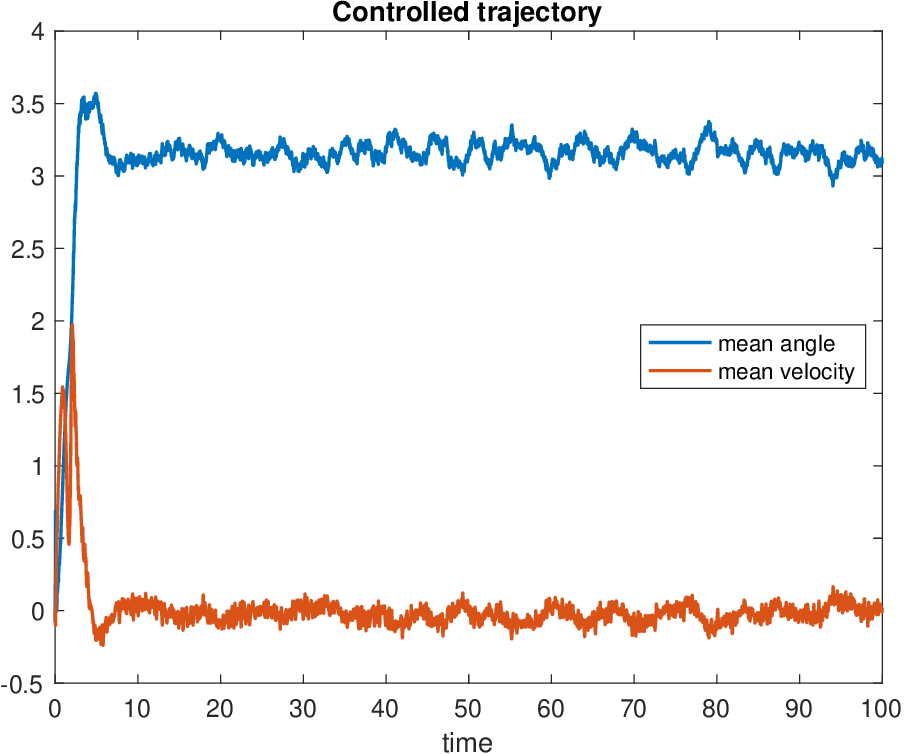}
$\quad$\includegraphics[width=0.47\textwidth]{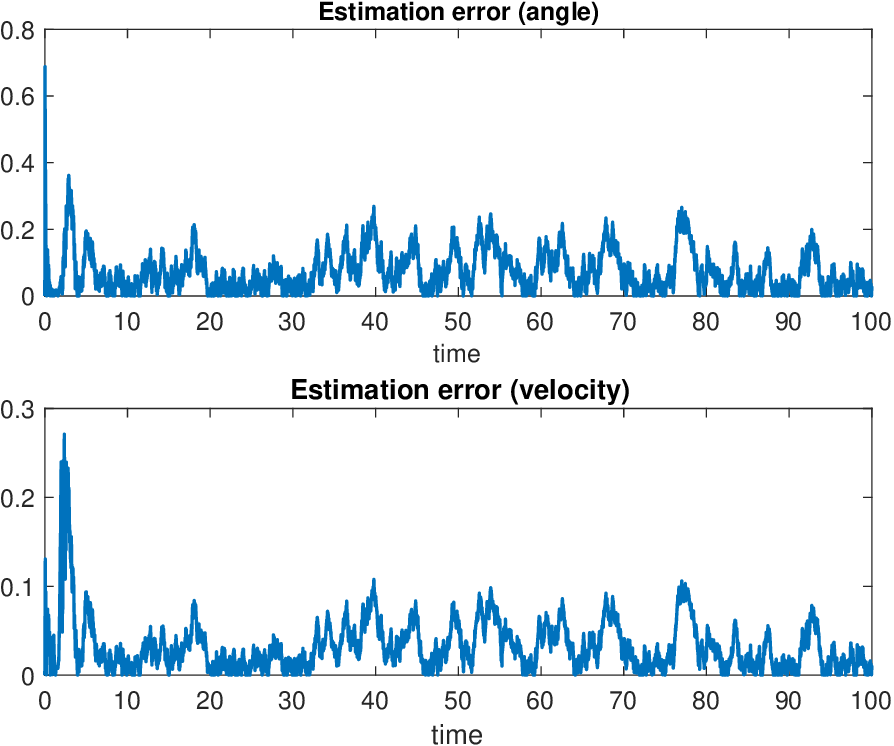}} 
\caption{Inverted pendulum.~Left panel: Controlled angle, $\theta_t$, and velocity, $v_t$, as a function of time. Starting from the stable equilibrium $(0,0)$, the control quickly stabilises the upper position (unstable) equilibrium $(\pi,0)$. Right panel: Estimation error in angle and velocity as a function of time.
} 
\label{fig7}
\end{figure}

We consider a controlled inverted pendulum with friction \cite{Meyn}. The state variable $x = (\theta,v)^{\rm T}$ is two-dimensional with equations of motion
\begin{subequations} \label{eq:pendulum}
\begin{align}
\dot{\theta}_t &= v_t,\\
\dot{v}_t &= -\sin (\theta_t) -\sigma v_t + \cos(\theta_t)U_t 
\end{align}
\end{subequations}
and $\sigma = 2$. We observe the angles $\theta_t^\ast$ of the physical twin continuously in time with measurement error variance $R = 0.01$ in (\ref{eq:cobs}).

We consider the running cost
\begin{equation} \label{eq:running cost IP}
c(x) = \frac{\rho}{2} \left( (\theta-\pi)^2
+ v^2 \right)
\end{equation}
with $\rho = 500$ and use the discount factor $\gamma = 1$
in (\ref{eq:cost}). Note that $(\pi,0)^{\rm T}$ is an unstable equilibrium point under $U_t \equiv 0$.

We again add diffusion (viscosity) to the digital twin setting $\Sigma = 0.1I$ in (\ref{eq:Hamiltonian density}). We initialize the particles at the stable equilibrium point $x_{\rm s} = (0,0)^{\rm T}$ with added Gaussian noise of variance $0.1I$; {\it i.e.}~$X_0^{(i)} \sim {\rm N}(x_{\rm s},0.1I)$. The initial momenta are  $P_0^{(i)} = 0$. The physical twin is initialized at $X_0^\dagger = x_{\rm s}$. The parameter $\delta$ in (\ref{eq:NWKR}) is set to $\delta = 0.1$ and the inflation parameter in (\ref{eq:inflation}) to $\sigma = 0.002$. The evolution equations (\ref{eq:interacting particle formulation}) have been simulated in 100,000 Euler time-steps with time-step $\Delta t = 0.001$ using $M=3$ particles.

The results can be found in Figure \ref{fig7}. It can be clearly seen that the twins leave the stable equilibrium and rapidly equilibrate at the unstable equilibrium even when using only $M=3$ particles in (\ref{eq:inflation}). We also display the resulting estimation errors in the angle and velocity. 

%%%%%%%%%%%%%%%%%%%%%%%%%%%%
%
\subsection{Controlled Lorenz-96 system} \label{sec:num_L96}
%
%%%%%%%%%%%%%%%%%%%%%%%%%%%%%

We now consider the Lorenz-96 system \cite{lorenz96} with $d_x = 40$ and the forcing set to $F=10$. This setting leads to rapidly changing and highly chaotic solutions. We observe every second component of $X_t \in \mathbb{R}^{d_x}$ continuously in time with measurement error $R =1$. The control task is to stabilize $x=0$ with $G=I$. We introduce the running cost
\begin{equation}
c(x) = \frac{\rho}{2} \|x\|^2
\end{equation}
with $\rho = 100$ and employ the discount factor $\gamma = 1$. The McKean--Pontryagin formulation (\ref{eq:Hamiltonian ODE continuous}) is implemented with $\Sigma = 0$ and, based on the approximation (\ref{eq:linear approximation}),
\begin{equation}
D_x \psi_t(X_t)\dot{X}_t \approx C_t^{px}\dot{X}_t.
\end{equation}
The ensemble size is set to $M=10$. Thus, $M\ll d_x$ and covariance localization and inflation are necessary. Covariance localization is based on the Gaspari--Cohn kernel \cite{gaspari99} with radius $r_{\rm loc} = 8$ \cite{reich2015probabilistic}. Localization is applied to both $C_{t}^{xx}$ and $C_{t}^{px}$. Given the step-size $\Delta t = 0.001$, the multiplicative inflation factor is set to $\gamma_{\rm inf} = 1+0.1\Delta t$; {\it i.e.},
\begin{equation}
    X_{t_n} \to m_{t_n} + \gamma_{\rm inf}(X_{t_n}-m_{t_n}), \qquad t_n = n\,\Delta t,
\end{equation}
before each data assimilation step with the ensemble Kalman--Bucy filter. 

\begin{figure}
\centerline{\includegraphics[width=0.47\textwidth]{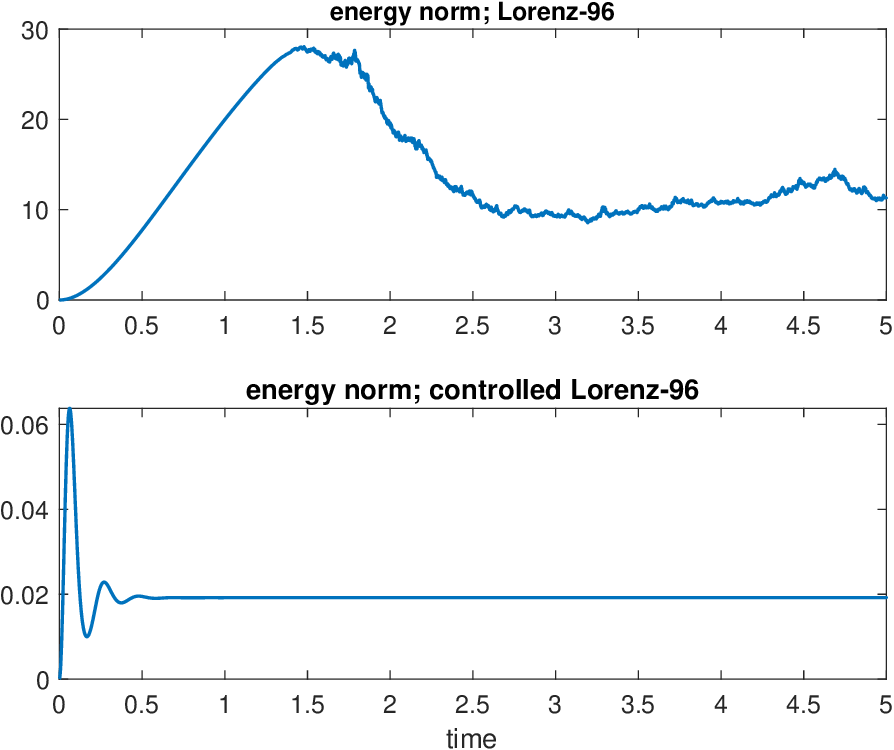}
} 
\caption{Lorenz-96.~Top panel: Time evolution of the energy norm (\ref{eq:energy norm}) for the partially observed Lorenz-96 system without control. Bottom panel: Same for the controlled system.
} 
\label{fig8}
\end{figure}

The system is initialized at $X_0 \sim {\rm N}(0,0.01I)$, $P_0 = 0$ and integrated up to time $T=5$. The time evolution of the scaled energy norm
\begin{equation} \label{eq:energy norm}
    e(t_n) = \frac{1}{2d_x} \|m^x_{t_n}\|^2, 
\end{equation}
$n = 0,\ldots,T/\Delta t$, can be found in Figure \ref{fig8} for both the uncontrolled and controlled Lorenz-96 system. The impact of control can be clearly seen and confirms that $x=0$ is effectively stabilized while the uncontrolled solutions develop into chaotic dynamics. 

%%%%%%%%%%%%%%%%%%%%%%%%%%%%%%%%%%%%
%
\section{Conclusions} \label{sec:conclusion}
%
%%%%%%%%%%%%%%%%%%%%%%%%%%%%%%%%%%%

Although the impact of digital twins on numerous application areas ranging from numerical weather prediction \citep{bauer2015,BSH21} to personalized medicine \citep{asch22,DT_medicine} and engineering \citep{DTWillcox,Willcox25} is beyond doubt, fundamental theoretical and algorithmic aspects still await their resolution. See the recent SIAM Report on the {\it Future of Computational Science} \citep{SIAM24} and National Academies of Sciences, Engineering and Medicine Report on {\it Foundational Research Gaps and Future Directions for Digital Twins} \citep{NA24}. 

In this paper, we have proposed an online digital twin formulation which combines the ensemble Kalman filter \cite{Evensenetal2022,CRS22} for data assimilation with the infinite horizon McKean--Pontryagin formulation for optimal control \cite{R25b}. We have thus addressed the problem of how to adapt a digital twin to incoming data from its physical twin while simultaneously providing control laws applicable to the physical twin. Although initial numerical studies for controlled Lorenz-63 and Lorenz-96 system as well as an inverted pendulum demonstrate the suitability of the proposed computational methodology, an in-depth mathematical investigation is left to future research. 

Our approach is related to model predictive control \cite{MPC} with the key difference that the filtering distribution and control laws are computed simultaneously and interactively. The proposed methodology thus also avoids an explicit application of the separation or equivalence principle \cite{Handel07,WW81}; while still providing an approximation to the partially observed Markov decision problem \cite{ASTROM1965,Bensoussan92}. See also Remark \ref{rem1}.

Extensions to high-dimensional evolution equations require additional approximations in (\ref{eq:interacting particle formulation}) such as covariance localization \cite{reich2015probabilistic,Evensenetal2022} and conditional independence in the approximation of the generator (\ref{eq:generator}) as introduced in \cite{gottwald2024localized}. Furthermore, one can restrict the class of 
functions $\beta_t(x)$ and $\psi_t(x)$ as, for example, in (\ref{eq:linear approximation}) \cite{R25b}.

\medskip

\paragraph{Acknowledgements}
This work has been funded by Deutsche Forschungsgemeinschaft (DFG) - Project-ID 318763901 - SFB1294. 

%%%%%%%%%%%%%%%%%%%%%%%%%%%%%%%%%

%\bibliographystyle{siamplain}
\bibliographystyle{elsarticle-num-names} 

%%%%%%%%%%%%%%%%%%%%%%%%%%%%%%%%%
%
\bibliography{bib-database.bib}

\begin{thebibliography}{52}
\providecommand{\natexlab}[1]{#1}
\providecommand{\url}[1]{\texttt{#1}}
\providecommand{\urlprefix}{URL }
\expandafter\ifx\csname urlstyle\endcsname\relax
  \providecommand{\doi}[1]{doi:\discretionary{}{}{}#1}\else
  \providecommand{\doi}[1]{doi:\discretionary{}{}{}\begingroup
  \urlstyle{rm}\url{#1}\endgroup}\fi
\providecommand{\bibinfo}[2]{#2}

\bibitem[{Meyn(2022)}]{Meyn}
\bibinfo{author}{S.~Meyn}, \bibinfo{title}{Control Systems and Reinforcement
  Learning}, \bibinfo{publisher}{Cambridge University Press},
  \bibinfo{address}{Cambridge}, \bibinfo{year}{2022}.

\bibitem[{Moerland et~al.(2023)Moerland, Broekens, Plaat, and
  Jonker}]{Model-based_reinforcement}
\bibinfo{author}{T.~M. Moerland}, \bibinfo{author}{J.~Broekens},
  \bibinfo{author}{A.~Plaat}, \bibinfo{author}{C.~M. Jonker},
  \bibinfo{title}{Model-based Reinforcement Learning: {A} Survey},
  \bibinfo{publisher}{Now Publishers}, \bibinfo{year}{2023}.

\bibitem[{Reich and Cotter(2015)}]{reich2015probabilistic}
\bibinfo{author}{S.~Reich}, \bibinfo{author}{C.~Cotter},
  \bibinfo{title}{Probabilistic Forecasting and {B}ayesian Data Assimilation},
  \bibinfo{publisher}{Cambridge University Press}, \bibinfo{year}{2015}.

\bibitem[{Law et~al.(2015)Law, Stuart, and Zygalakis}]{law2015data}
\bibinfo{author}{K.~Law}, \bibinfo{author}{A.~Stuart},
  \bibinfo{author}{K.~Zygalakis}, \bibinfo{title}{Data Assimilation},
  \bibinfo{journal}{Cham, Switzerland: Springer} \bibinfo{volume}{214}.

\bibitem[{\r{A}str\"om(1965)}]{ASTROM1965}
\bibinfo{author}{K.~J. \r{A}str\"om}, \bibinfo{title}{Optimal control of
  {M}arkov processes with incomplete state information},
  \bibinfo{journal}{Journal of Mathematical Analysis and Applications}
  \bibinfo{volume}{10} (\bibinfo{year}{1965}) \bibinfo{pages}{174--205}.

\bibitem[{Bensoussan(1992)}]{Bensoussan92}
\bibinfo{author}{A.~Bensoussan}, \bibinfo{title}{Stochastic Control of
  Partially Observable Systems}, \bibinfo{publisher}{Cambridge University
  Press}, \bibinfo{address}{Cambridge}, \bibinfo{year}{1992}.

\bibitem[{{van Handel}(2007)}]{Handel07}
\bibinfo{author}{R.~{van Handel}}, \bibinfo{title}{Stochastic Calculus,
  Filtering, and Stochastic Control}, \bibinfo{howpublished}{Lecture Notes,
  Caltech}, \bibinfo{year}{2007}.

\bibitem[{Nisio(2015)}]{Nisio15}
\bibinfo{author}{M.~Nisio}, \bibinfo{title}{Stochastic Control Theory},
  \bibinfo{publisher}{Springer}, \bibinfo{address}{Tokyo},
  \bibinfo{edition}{2nd} edn., \bibinfo{year}{2015}.

\bibitem[{Silver and Veness(2010)}]{Silver10}
\bibinfo{author}{D.~Silver}, \bibinfo{author}{J.~Veness},
  \bibinfo{title}{{Monte-Carlo Planning in Large POMDPs}}, in:
  \bibinfo{editor}{J.~Lafferty}, \bibinfo{editor}{C.~Williams},
  \bibinfo{editor}{J.~Shawe-Taylor}, \bibinfo{editor}{R.~Zemel},
  \bibinfo{editor}{A.~Culotta} (Eds.), \bibinfo{booktitle}{Advances in Neural
  Information Processing Systems}, vol.~\bibinfo{volume}{23},
  \bibinfo{publisher}{Curran Associates, Inc.}, \bibinfo{pages}{2164--2172},
  \bibinfo{year}{2010}.

\bibitem[{{van de Water} and Willems(1981)}]{WW81}
\bibinfo{author}{H.~{van de Water}}, \bibinfo{author}{J.~C. Willems},
  \bibinfo{title}{The certainty equivalence property in stochastic control
  theory}, \bibinfo{journal}{IEEE Trans. Autom. Contr.} \bibinfo{volume}{AC-26}
  (\bibinfo{year}{1981}) \bibinfo{pages}{1080--1087}.

\bibitem[{Littman et~al.(1995)Littman, Cassandra, and Kaelbling}]{QMDP}
\bibinfo{author}{M.~L. Littman}, \bibinfo{author}{A.~R. Cassandra},
  \bibinfo{author}{L.~P. Kaelbling}, \bibinfo{title}{Learning policies for
  partially observable environments: {S}caling up}, in:
  \bibinfo{booktitle}{International Conference on Machine Learning},
  \bibinfo{pages}{362--370}, \bibinfo{year}{1995}.

\bibitem[{Ross et~al.(2008)Ross, Pineau, Paquet, and Chaib-draa}]{RPPC08}
\bibinfo{author}{S.~Ross}, \bibinfo{author}{J.~Pineau},
  \bibinfo{author}{S.~Paquet}, \bibinfo{author}{B.~Chaib-draa},
  \bibinfo{title}{Online Planning algorithms for {POMDP}s},
  \bibinfo{journal}{Journal of Artificial Inelligence Research}
  \bibinfo{volume}{32} (\bibinfo{year}{2008}) \bibinfo{pages}{663--704}.

\bibitem[{Zheng et~al.(2022)Zheng, Ridderhof, Tsiotras, and
  Agha-Mohammadi}]{9811560}
\bibinfo{author}{D.~Zheng}, \bibinfo{author}{J.~Ridderhof},
  \bibinfo{author}{P.~Tsiotras}, \bibinfo{author}{A.-A. Agha-Mohammadi},
  \bibinfo{title}{Belief Space Planning: {A} Covariance Steering Approach}, in:
  \bibinfo{booktitle}{2022 International Conference on Robotics and Automation
  (ICRA)}, \bibinfo{pages}{11051--11057}, \bibinfo{year}{2022}.

\bibitem[{Zheng et~al.(2024)Zheng, Ridderhof, Zhang, Tsiotras, and
  Agha-Mohammadi}]{BSP24}
\bibinfo{author}{D.~Zheng}, \bibinfo{author}{J.~Ridderhof},
  \bibinfo{author}{Z.~Zhang}, \bibinfo{author}{P.~Tsiotras},
  \bibinfo{author}{A.-A. Agha-Mohammadi}, \bibinfo{title}{{CS-BRM}: {A}
  Probabilistic RoadMap for Consistent Belief Space Planning With Reachability
  Guarantees}, \bibinfo{journal}{IEEE Transactions on Robotics}
  (\bibinfo{year}{2024}) \bibinfo{pages}{1630--1649}.

\bibitem[{Kappen(2005)}]{Kappen05}
\bibinfo{author}{H.~Kappen}, \bibinfo{title}{Path integrals and symmetry
  breaking for optimal control theory}, \bibinfo{journal}{J Statistical
  Mechanics: Theory and Experiments} \bibinfo{volume}{11}
  (\bibinfo{year}{2005}) \bibinfo{pages}{11011}.

\bibitem[{Kappen et~al.(2012)Kappen, G\'omez, and Opper}]{KVO12}
\bibinfo{author}{H.~J. Kappen}, \bibinfo{author}{V.~G\'omez},
  \bibinfo{author}{M.~Opper}, \bibinfo{title}{Optimal control as a graphical
  model inference problem}, \bibinfo{journal}{Machine Learning}
  \bibinfo{volume}{87} (\bibinfo{year}{2012}) \bibinfo{pages}{159--182}.

\bibitem[{Abdulsamad et~al.(2025)Abdulsamad, Iqbal, and S\"arkk\"a}]{AIS25}
\bibinfo{author}{H.~Abdulsamad}, \bibinfo{author}{S.~Iqbal},
  \bibinfo{author}{S.~S\"arkk\"a}, \bibinfo{title}{Sequential {M}onte {C}arlo
  for policy optimization in continuous {POMDP}s}, \bibinfo{type}{Tech. Rep.},
  \bibinfo{institution}{arXiv:2505.16732}, \bibinfo{year}{2025}.

\bibitem[{Kawasaki and Kotsuki(2024)}]{KK24}
\bibinfo{author}{F.~Kawasaki}, \bibinfo{author}{S.~Kotsuki},
  \bibinfo{title}{Leading the {L}orenz 63 system towards the prescribed regime
  by model predictive control coupled with data assimilation},
  \bibinfo{journal}{Nonlin. Processes Geophys.} \bibinfo{volume}{31}
  (\bibinfo{year}{2024}) \bibinfo{pages}{319--333}.

\bibitem[{Reich(2025)}]{R25a}
\bibinfo{author}{S.~Reich}, \bibinfo{title}{{Ensemble {K}alman--{B}ucy
  filtering for nonlinear model predictive control}}, \bibinfo{type}{Tech.
  Rep.}, \bibinfo{institution}{arXiv:2503.12474}, \bibinfo{year}{2025}.

\bibitem[{Evensen et~al.(2022)Evensen, Vossepoel, and {van
  Leeuwen}}]{Evensenetal2022}
\bibinfo{author}{G.~Evensen}, \bibinfo{author}{F.~C. Vossepoel},
  \bibinfo{author}{P.~J. {van Leeuwen}}, \bibinfo{title}{Data Assimilation
  Fundamentals: {A} unified Formulation of the State and Parameter Estimation
  Problem}, \bibinfo{publisher}{Springer Nature Switzerland AG},
  \bibinfo{address}{Cham, Switzerland}, \bibinfo{year}{2022}.

\bibitem[{Opper and Reich(2025)}]{R25b}
\bibinfo{author}{M.~Opper}, \bibinfo{author}{S.~Reich}, \bibinfo{title}{On a
  mean-field {P}ontryagin minimum principle for stochastic optimal control},
  \bibinfo{type}{Tech. Rep.}, \bibinfo{institution}{arXiv:2506.10506},
  \bibinfo{year}{2025}.

\bibitem[{Calvello et~al.(2025)Calvello, Reich, and Stuart}]{CRS22}
\bibinfo{author}{E.~Calvello}, \bibinfo{author}{S.~Reich},
  \bibinfo{author}{A.~M. Stuart}, \bibinfo{title}{Ensemble {K}alman methods:
  {A} mean field perspective}, \bibinfo{journal}{Acta Numerica}
  \bibinfo{volume}{34} (\bibinfo{year}{2025}) \bibinfo{pages}{123--291}.

\bibitem[{Lorenz(1996)}]{lorenz96}
\bibinfo{author}{E.~Lorenz}, \bibinfo{title}{Predictibility: {A} problem partly
  solved}, in: \bibinfo{booktitle}{Proc. Seminar on Predictibility},
  vol.~\bibinfo{volume}{1}, \bibinfo{address}{ECMWF, Reading, Berkshire, UK},
  \bibinfo{pages}{1--18}, \bibinfo{year}{1996}.

\bibitem[{Bergemann and Reich(2012)}]{bergemann2012ensemble}
\bibinfo{author}{K.~Bergemann}, \bibinfo{author}{S.~Reich}, \bibinfo{title}{An
  ensemble {K}alman--{B}ucy filter for continuous data assimilation},
  \bibinfo{journal}{Meteorologische Zeitschrift}
  \bibinfo{volume}{21}~(\bibinfo{number}{3}) (\bibinfo{year}{2012})
  \bibinfo{pages}{213}.

\bibitem[{Coghi et~al.(2023)Coghi, Nilssen, N\"usken, and Reich}]{CNN2021}
\bibinfo{author}{M.~Coghi}, \bibinfo{author}{T.~Nilssen},
  \bibinfo{author}{N.~N\"usken}, \bibinfo{author}{S.~Reich},
  \bibinfo{title}{Rough {M}c{K}ean--{V}lasov dynamics for robust ensemble
  {K}alman filtering}, \bibinfo{journal}{Ann. Appl. Probab.}
  \bibinfo{volume}{33 (6B)} (\bibinfo{year}{2023}) \bibinfo{pages}{5693--5752}.

\bibitem[{Chopin and Papaspiliopoulos(2020)}]{chopin:20}
\bibinfo{author}{N.~Chopin}, \bibinfo{author}{O.~Papaspiliopoulos},
  \bibinfo{title}{An Introduction to Sequential {M}onte {C}arlo},
  \bibinfo{publisher}{Springer Nature Switzerland AG}, \bibinfo{address}{Cham,
  Switzerland}, \bibinfo{year}{2020}.

\bibitem[{Snyder et~al.(2008)Snyder, Bengtsson, Bickel, and Anderson}]{SBBA08}
\bibinfo{author}{C.~Snyder}, \bibinfo{author}{T.~Bengtsson},
  \bibinfo{author}{P.~Bickel}, \bibinfo{author}{J.~Anderson},
  \bibinfo{title}{Obstacles to High-Dimensional Particle Filtering},
  \bibinfo{journal}{Monthly Weather Review} \bibinfo{volume}{136}
  (\bibinfo{year}{2008}) \bibinfo{pages}{4629--4640}.

\bibitem[{Bergemann and Reich(2010)}]{bergemann2010mollified}
\bibinfo{author}{K.~Bergemann}, \bibinfo{author}{S.~Reich}, \bibinfo{title}{A
  mollified ensemble {K}alman filter}, \bibinfo{journal}{Quarterly Journal of
  the Royal Meteorological Society}
  \bibinfo{volume}{136}~(\bibinfo{number}{651}) (\bibinfo{year}{2010})
  \bibinfo{pages}{1636--1643}.

\bibitem[{Carmona(2016)}]{Carmona}
\bibinfo{author}{R.~Carmona}, \bibinfo{title}{Lectures on {BSDE}s, Stochastic
  Control, and Stochastic Differential Games with Financial Applications},
  \bibinfo{publisher}{SIAM}, \bibinfo{address}{Philadelphia},
  \bibinfo{year}{2016}.

\bibitem[{Pavliotis(2016)}]{Pavliotis2016}
\bibinfo{author}{G.~A. Pavliotis}, \bibinfo{title}{Stochastic Processes and
  Applications}, \bibinfo{publisher}{Springer Verlag}, \bibinfo{address}{New
  York}, \bibinfo{year}{2016}.

\bibitem[{Yang et~al.(2013)Yang, Mehta, and Meyn}]{meyn13}
\bibinfo{author}{T.~Yang}, \bibinfo{author}{P.~G. Mehta},
  \bibinfo{author}{S.~P. Meyn}, \bibinfo{title}{Feedback particle filter},
  \bibinfo{journal}{IEEE Trans. Automat. Control}
  \bibinfo{volume}{58}~(\bibinfo{number}{10}) (\bibinfo{year}{2013})
  \bibinfo{pages}{2465--2480}, ISSN \bibinfo{issn}{0018-9286}.

\bibitem[{Pointryagin et~al.(1962)Pointryagin, Boltyanskii, Gamkrelidze, and
  Mihchenko}]{pontryagin}
\bibinfo{author}{L.~Pointryagin}, \bibinfo{author}{V.~Boltyanskii},
  \bibinfo{author}{R.~Gamkrelidze}, \bibinfo{author}{E.~Mihchenko},
  \bibinfo{title}{The Mathematical Theory of Optimal Processes},
  \bibinfo{publisher}{John Wiley \& Sons}, \bibinfo{address}{New York},
  \bibinfo{year}{1962}.

\bibitem[{Bensoussan(2018)}]{Bensoussan}
\bibinfo{author}{A.~Bensoussan}, \bibinfo{title}{Estimation and Control of
  Dynamical Systems}, \bibinfo{publisher}{Springer}, \bibinfo{address}{Cham},
  \bibinfo{year}{2018}.

\bibitem[{Mehta and Meyn(2013)}]{PM13}
\bibinfo{author}{P.~G. Mehta}, \bibinfo{author}{S.~P. Meyn}, \bibinfo{title}{A
  feedback particle filter-based approach to optimal control with partial
  observations}, in: \bibinfo{booktitle}{52nd IEEE Conference on Decision and
  Control}, \bibinfo{pages}{3121--3127}, \bibinfo{year}{2013}.

\bibitem[{Bierens(1994)}]{Bierens}
\bibinfo{author}{H.~J. Bierens}, \bibinfo{title}{The {N}adaraya--{W}atson
  kernel regression function estimator}, in: \bibinfo{booktitle}{Topics in
  Advanced Econometrics}, \bibinfo{publisher}{Cambridge University Press},
  \bibinfo{address}{New York}, \bibinfo{pages}{212--247}, \bibinfo{year}{1994}.

\bibitem[{Wormell and Reich(2021)}]{WR20}
\bibinfo{author}{C.~L. Wormell}, \bibinfo{author}{S.~Reich},
  \bibinfo{title}{Spectral convergence of diffusion maps: {I}mproved error
  bounds and an alternative normalisation}, \bibinfo{journal}{SIAM J. Numer.
  Anal.} \bibinfo{volume}{59} (\bibinfo{year}{2021})
  \bibinfo{pages}{1687--1734}.

\bibitem[{Gottwald et~al.(2025)Gottwald, Li, Reich, and Marzouk}]{GLRY24}
\bibinfo{author}{G.~A. Gottwald}, \bibinfo{author}{F.~Li},
  \bibinfo{author}{S.~Reich}, \bibinfo{author}{Y.~Marzouk},
  \bibinfo{title}{{Stable generative modeling using Schr\"odinger bridges}},
  \bibinfo{journal}{Phil. Trans. R. Soc. A} \bibinfo{volume}{383}
  (\bibinfo{year}{2025}) \bibinfo{pages}{20240332}.

\bibitem[{Marshall and Coifman(2019)}]{MarshallCoifman}
\bibinfo{author}{N.~F. Marshall}, \bibinfo{author}{R.~R. Coifman},
  \bibinfo{title}{Manifold learning with bi-stochastic kernels},
  \bibinfo{journal}{IMA J. Appl. Maths.} \bibinfo{volume}{84}
  (\bibinfo{year}{2019}) \bibinfo{pages}{455--482}.

\bibitem[{Gottwald and Reich(2026)}]{gottwald2024localized}
\bibinfo{author}{G.~A. Gottwald}, \bibinfo{author}{S.~Reich},
  \bibinfo{title}{Localized {S}chr\"odinger Bridge Sampler},
  \bibinfo{journal}{J.~Comput.~Phys} \bibinfo{volume}{548}
  (\bibinfo{year}{2026}) \bibinfo{pages}{114583}.

\bibitem[{Lorenz(1963)}]{lorenz1963deterministic}
\bibinfo{author}{E.~N. Lorenz}, \bibinfo{title}{Deterministic nonperiodic
  flow}, \bibinfo{journal}{Journal of the Atmospheric Sciences}
  \bibinfo{volume}{20}~(\bibinfo{number}{2}) (\bibinfo{year}{1963})
  \bibinfo{pages}{130--141}.

\bibitem[{Findeisen and Allg\"ower(2002)}]{Allgower}
\bibinfo{author}{R.~Findeisen}, \bibinfo{author}{F.~Allg\"ower},
  \bibinfo{title}{An Introduction to Nonlinear Model Predictive Control}, in:
  \bibinfo{editor}{A.~G. Jager, de}, \bibinfo{editor}{H.~J. Zwart} (Eds.),
  \bibinfo{booktitle}{Systems and control : 21th Benelux meeting 2002},
  \bibinfo{publisher}{Technische Universiteit Eindhoven},
  \bibinfo{pages}{119--141}, \bibinfo{year}{2002}.

\bibitem[{Rawlings et~al.(2018)Rawlings, Mayne, and Diehl}]{MPC}
\bibinfo{author}{J.~B. Rawlings}, \bibinfo{author}{D.~Q. Mayne},
  \bibinfo{author}{M.~M. Diehl}, \bibinfo{title}{Model Predictive Control:
  {T}heory, Computation, and Design}, \bibinfo{publisher}{Nob Hill Publishing},
  \bibinfo{address}{Madison}, \bibinfo{edition}{2nd} edn.,
  \bibinfo{year}{2018}.

\bibitem[{Crandall and Lions(1983)}]{CL83}
\bibinfo{author}{M.~G. Crandall}, \bibinfo{author}{P.-L. Lions},
  \bibinfo{title}{VISCOSITY SOLUTIONS OF {H}AMILTON-{J}ACOBI EQUATIONS},
  \bibinfo{journal}{Transactions of the American Mathematical Society}
  \bibinfo{volume}{277}.

\bibitem[{Gaspari and Cohn(1999)}]{gaspari99}
\bibinfo{author}{G.~Gaspari}, \bibinfo{author}{S.~Cohn},
  \bibinfo{title}{Construction of correlation functions in two and three
  dimensions}, \bibinfo{journal}{Q. J. Royal Meteorological Soc.}
  \bibinfo{volume}{125} (\bibinfo{year}{1999}) \bibinfo{pages}{723--757}.

\bibitem[{Bauer et~al.(2015)Bauer, Thorpe, and Brunet}]{bauer2015}
\bibinfo{author}{P.~Bauer}, \bibinfo{author}{A.~Thorpe},
  \bibinfo{author}{G.~Brunet}, \bibinfo{title}{The quiet revolution of
  numerical weather prediction}, \bibinfo{journal}{Nature}
  \bibinfo{volume}{525} (\bibinfo{year}{2015}) \bibinfo{pages}{47--55}.

\bibitem[{Bauer et~al.(2021)Bauer, Stevens, and Hazeleger}]{BSH21}
\bibinfo{author}{P.~Bauer}, \bibinfo{author}{B.~Stevens},
  \bibinfo{author}{W.~Hazeleger}, \bibinfo{title}{A digital twin of {E}arth for
  the green transition}, \bibinfo{journal}{Nature Climate Change}
  \bibinfo{volume}{11} (\bibinfo{year}{2021}) \bibinfo{pages}{80--83}.

\bibitem[{Asch(2022)}]{asch22}
\bibinfo{author}{M.~Asch}, \bibinfo{title}{A Toolbox for Digital Twins: {F}rom
  Model-Based to Data-Driven}, \bibinfo{publisher}{SIAM},
  \bibinfo{address}{Philadelphia}, \bibinfo{year}{2022}.

\bibitem[{Sel et~al.(2025)Sel, Hawkins-Daarud, Chaudhuri, Osman, Bahai,
  Paydarfar, Willcox, Chung, and Jafari}]{DT_medicine}
\bibinfo{author}{K.~Sel}, \bibinfo{author}{A.~Hawkins-Daarud},
  \bibinfo{author}{A.~Chaudhuri}, \bibinfo{author}{D.~Osman},
  \bibinfo{author}{A.~Bahai}, \bibinfo{author}{D.~Paydarfar},
  \bibinfo{author}{K.~Willcox}, \bibinfo{author}{C.~Chung},
  \bibinfo{author}{R.~Jafari}, \bibinfo{title}{Survey and perspective on
  verification, validation, and uncertainty quantification of digital twins for
  precision medicine}, \bibinfo{journal}{npj Digital Medicine}
  \bibinfo{volume}{8} (\bibinfo{year}{2025}) \bibinfo{pages}{40}.

\bibitem[{Niederer et~al.(2021)Niederer, Sacks, Girolami, and
  Willcox}]{DTWillcox}
\bibinfo{author}{S.~A. Niederer}, \bibinfo{author}{M.~S. Sacks},
  \bibinfo{author}{M.~Girolami}, \bibinfo{author}{K.~Willcox},
  \bibinfo{title}{Scaling digital twins from the artisanal to the industrial},
  \bibinfo{journal}{Nature Computational Science} \bibinfo{volume}{1}
  (\bibinfo{year}{2021}) \bibinfo{pages}{313--320}.

\bibitem[{Henao-Garcia et~al.(2025)Henao-Garcia, Kapteyn, Willcox, Tezzele,
  Castroviejo-Fernandez, Kim, Ambrosino, Kolmanovsky, Basu, Jirwankar, and
  Sanfelice}]{Willcox25}
\bibinfo{author}{S.~Henao-Garcia}, \bibinfo{author}{M.~Kapteyn},
  \bibinfo{author}{K.~E. Willcox}, \bibinfo{author}{M.~Tezzele},
  \bibinfo{author}{M.~Castroviejo-Fernandez}, \bibinfo{author}{T.~Kim},
  \bibinfo{author}{M.~Ambrosino}, \bibinfo{author}{I.~Kolmanovsky},
  \bibinfo{author}{H.~Basu}, \bibinfo{author}{P.~Jirwankar},
  \bibinfo{author}{R.~Sanfelice}, \bibinfo{title}{Digital-Twin-Enabled
  Multi-Spacecraft On-Orbit Operations}, AIAA SCITECH 2025 Forum,
  \bibinfo{pages}{1432}, \doi{\bibinfo{doi}{10.2514/6.2025-1432}},
  \bibinfo{year}{2025}.

\bibitem[{SIAM(2024)}]{SIAM24}
\bibinfo{author}{SIAM}, \bibinfo{title}{{Report on the Future of Computational
  Science}}, \bibinfo{type}{Tech. Rep.}, \bibinfo{institution}{SIAM
  Publishing}, \bibinfo{address}{Philadelphia}, \bibinfo{year}{2024}.

\bibitem[{{National Academy of Sciences}(2024)}]{NA24}
\bibinfo{author}{{National Academy of Sciences}}, \bibinfo{title}{Foundational
  Research Gaps and Future Directions for Digital Twins}, \bibinfo{type}{Tech.
  Rep.}, \bibinfo{institution}{Washington, DC: The National Academies Press},
  \bibinfo{year}{2024}.

\end{thebibliography}
%

%%%%%%%%%%%%%%%%%%%%%%%%%%%%

\end{document}